\documentclass[11pt]{article}
\usepackage{amsfonts}

\usepackage{graphics}
\usepackage{indentfirst}
\usepackage{cite}
\usepackage{latexsym}
\usepackage{amsmath,amsthm}
\usepackage{amssymb}
\usepackage[dvips]{epsfig}
\usepackage{amscd}
\usepackage{mathrsfs}

\usepackage{color}
\newtheorem{theorem}{Theorem}[section]
\newtheorem{remark}{Remark}[section]

\newtheorem{definition}{Definition}[section]
\newtheorem{lemma}[theorem]{Lemma}

\newtheorem{proposition}[theorem]{Proposition}

\newcommand{\n}{\rho}

\newcommand{\lm}{\lambda}

\renewcommand{\div}{ {\rm div }  }

\newcommand{\na}{\nabla }

\newcommand{\pa}{\partial}

\newcommand{\bt}{\begin{theorem}}
\newcommand{\bl}{\begin{lemma}}
\newcommand{\el}{\end{lemma}}
\newcommand{\et}{\end{theorem}}
\newcommand{\ga}{\gamma}

\newcommand{\OM}{\Omega}

\newcommand{\curl}{{\rm curl} }

\newcommand{\de}{\delta}

\newcommand{\la}{\label}
\newcommand{\si}{\sigma}

\newcommand{\om}{\Omega}
\newcommand{\ol}{\overline}

\newcommand{\bn}{\begin{eqnarray}}
\newcommand{\en}{\end{eqnarray}}
\newcommand{\bnn}{\begin{eqnarray*}}
\newcommand{\enn}{\end{eqnarray*}}

\newcommand{\bnnn}{\begin{eqnarray*}}
\newcommand{\ennn}{\end{eqnarray*}}
\newcommand{\ben}{\begin{enumerate}}
\newcommand{\een}{\end{enumerate}}

\newcommand{\rs}{\rho^*}

\newcommand{\T}{\mathbb{T}}

\newcommand{\du}{\dot{u}}

\newcommand{\ba}{\begin{aligned}}
\newcommand{\ea}{\end{aligned}}
\newcommand{\be}{\begin{equation}}
\newcommand{\ee}{\end{equation}}

\def\p{\partial}
\def\norm[#1]#2{\|#2\|_{#1}}

\def\lam{\lambda}
\def\ep{\varepsilon}

\def\o{\omega}
\def\rr{\mathbb{R}^2}

\makeatletter      
\@addtoreset{equation}{section}
\makeatother       

\title{Global Existence and Incompressible Limit for Compressible Navier-Stokes Equations in Bounded Domains with Large Bulk Viscosity and Large Initial Data}

\date{}

\author{$\text{Qinghao L{\small EI}}^{a,b}, \text{Chengfeng X{\small{IONG}}}^{a,b}\thanks{Email addresses:  leiqinghao22@mails.ucas.ac.cn (Q. H. Lei), xiongchengfeng20@mails.ucas.ac.cn (C. F. Xiong) }$\\
	a. School of Mathematical Sciences,\\ University of Chinese Academy of Sciences,
	Beijing 100049, P. R. China;\\
	b. Institute of Applied Mathematics,\\ Academy of Mathematics and Systems Science,\\
	Chinese Academy of Sciences, Beijing 100190, P. R. China.}

\begin{document}
\maketitle

\begin{abstract}
We investigate the barotropic compressible Navier-Stokes equations with Navier-slip boundary conditions in a general two-dimensional bounded simply connected domain.
For initial data allowing vacuum, we establish the global existence and exponential decay of weak, strong, and classical solutions when the bulk viscosity coefficient is sufficiently large, without any restrictions on the size of the initial data.
Furthermore, we prove that, as the bulk viscosity coefficient tends to infinity, solutions of the compressible Navier-Stokes equations converge to those of the inhomogeneous incompressible Navier-Stokes equations. 
In particular, we establish the incompressible limit of weak solutions without assuming that the initial velocity field is divergence-free. 
\end{abstract}

\section{Introduction and main results}
We study the two-dimensional barotropic compressible
Navier-Stokes equations which read as follows:
\be\ba\la{ns}
\begin{cases}
  \rho_t + \div(\rho u) = 0,\\
  (\n u)_t + \div(\n u\otimes u) -\mu \Delta u - (\mu + \lm)\na\div u
    +\na P = 0,
\end{cases}
\ea\ee
where $t \ge 0$ is time, $x \in \OM \subset \rr$ is the spatial coordinate,
$\n=\n(x,t)$ and $u(x,t)=(u^1(x,t),u^2(x,t))$ represent the 
density and velocity of the compressible flow respectively,
and the pressure $P$ is given by
\be\la{i1}
  P=a\n^\ga,
\ee
with consants $a>0,{\ga > 1}$. Without loss of generality, it is assumed that
$a=1$. The shear viscosity coefficient $\mu$ and bulk viscosity coefficient $\lam$ satisfy
the physical restrictions:
\be\la{i2}
\mu>0,\quad \mu+\lam\geq 0.
\ee
For later purpose, we set 
\be\la{nu}
\nu \triangleq 2\mu + \lam,
\ee
which together with (\ref{i2}) yields that 
\be\la{numu}\nu\geq\mu.\ee
In this paper, we assume that $\OM$ is a simply connected bounded $C^{2,1}$ domain in $\rr$.
In addition, the system is subject to the given initial data
\be\la{i30}
\n(x,0)=\n_0(x),\quad \n u(x,0)=\n_0u_0(x), \quad x\in \OM,
\ee
and Navier-slip boundary conditions:
\be\la{i3}\ba
u \cdot n = 0, \,\,\curl u=-Au \cdot n^\bot \,\,\,\text{on}\,\,\, \partial\Omega,
\ea\ee
where $A$ is a non-negative smooth function on the boundary, $n=(n_1,n_2)$ denotes the unit outer normal vector of the boundary $\partial \Omega$, while $n^\bot$ is the unit tangential vector on $\partial \Omega$ denoted by
\be\ba\nonumber
n^\bot\triangleq (-n_2,n_1).
\ea\ee
In addition, we carry out a smooth extension of $n$, $n^\bot$, and $A$ to $\overline\Omega$.
Considering that there are various ways to make this extension, 
we select a single approach and stick with it during the {rest} of the paper.

It is obvious that the total mass of smooth enough solutions
of (\ref{ns}) is conserved through the evolution, that is, for all
$t>0$,
\be\ba\la{conv}
\int_{\OM} \n dx = \int_{\OM}\n_0 dx.
\ea\ee
Without loss of generality, we shall assume that 
\be\ba\la{asm}
\int_{\OM} \n_0 dx = 1.
\ea\ee

There is a lot of literature concerning the global existence of weak and classical solutions to (\ref{ns}).
For the one-dimensional case, numerous researchers have derived
comprehensive results (see \cite{H4,KS,S1,S2} and references therein).
For the multi-dimensional case, Nash \cite{N} and Serrin \cite{S} established the
local existence and uniqueness of classical solutions respectively with the absence of vacuum. 
Furthermore, when the initial density is not required to be strictly positive 
and can vanish on open sets, the local existence
and uniqueness of strong solutions were demonstrated in \cite{CCK,CK,CK2,SS,LLL} 
and references therein.
The first result of global classical solutions was obtained by
Matsumura--Nishida\cite{MN1}, in which the initial data were required to be close
to a non-vacuum equilibrium in some Sobolev space $H^s$.
Subsequently, Hoff \cite{H1,H2,H3} investigated the problem for discontinuous
initial data and developed new a priori estimates for the material derivative $\dot{u}$.
Regarding the global existence of weak solutions for arbitrarily large initial data,
the major breakthrough was achieved by Lions \cite{L2}.
Under the finite initial energy assumption,
he established the global existence of weak solutions with vacuum 
when $\ga$ is suitably large. Specifically, 
for three-dimensional cases, $\ga \ge \frac{9}{5}$. 
These results were later refined by Feireisl-Novotn\'y-Petzeltov\'a 
\cite{FNP} to $\ga > \frac{3}{2}$.
Recently, Huang-Li-Xin\cite{HLX2} and Li-Xin\cite{LX2} proved the global
existence and uniqueness of classical solutions to the Cauchy problem
in three-dimensional and two-dimensional spaces respectively.
Their results require the initial energy to be sufficiently small
while allowing for large oscillations of the density,
and the initial density may contain vacuum and even have compact support.
Later, Cai-Li \cite{CL} generalized the above results to bounded domains 
with the velocity field subject to Navier-slip boundary conditions.

More recently, Danchin-Mucha \cite{DM} established the global existence of weak solutions
when the bulk viscosity is sufficiently large and $\nu^{1/2} \|\div u_0\|_{L^2}$ is bounded.
Furthermore, they also demonstrated that as the bulk viscosity tends to infinity,
the weak solutions converge to solutions of the
inhomogeneous incompressible Navier-Stokes equations.
Very recently, in our work \cite{LX3}, we removed the restriction on $\nu^{1/2} \|\div u_0\|_{L^2}$ of \cite{DM},
establishing global existence solely under the assumption that the bulk viscosity coefficient is suitably large.
However, these results were established only for the periodic case.
The main aim of this paper is to establish the global existence, exponential decay, and incompressible limit of weak, strong, and classical solutions for the two-dimensional compressible Navier-Stokes equations in general bounded simply connected domains.
These results are obtained without any restrictions on the size of the initial data, provided that the bulk viscosity coefficient is sufficiently large.

Before stating the main results, we first explain the notations
and conventions used throughout this paper. We denote
\be\ba\nonumber
\int f dx = \int_{\OM} fdx,\quad
\ol{f}=\frac{1}{|\OM|}\int f dx.
\ea\ee
For a positive integer $s$ and $1\leq r\leq\infty$, we also denote the standard Lebesgue and Sobolev
spaces as follows:
\be\ba\nonumber
\begin{cases}
L^r =L^r(\OM),\quad W^{s,r} =W^{s,r}(\OM),\quad H^s =W^{s,2}, \\
\tilde{H}^1 =\{u\in H^1(\Omega )\vert u \cdot n=0, \mathrm{curl}u =-Au \cdot n^\bot \text{ on } \partial \Omega \}.
\end{cases}
\ea\ee
Next, the material derivative and the transpose gradient are given by
\be\ba\nonumber
\frac{D}{Dt}f=\dot{f} \triangleq f_t + u\cdot\na f, \quad \na^{\bot} \triangleq (-\pa_2,\pa_1).
\ea\ee
The initial total energy is defined as follows:
\be\la{e0}\ba
E_0 \triangleq \int \left( \frac{1}{2} \rho_0 |u_0|^2 + \frac{1}{\ga-1}\n^\ga_0 \right) dx.
\ea\ee
Additionally, we define
\be\la{xd}\ba
\o \triangleq \na^\bot \cdot u=\pa_1 u^2 - \pa_2 u^1.
\ea\ee

Then we provide the definition of weak and strong solutions to (\ref{ns}).
\begin{definition}
If $(\n,u)$ satisfies \eqref{ns} in the sense of distribution, then we call $(\n,u)$ a weak solution.
In addition, for a weak solution if
all derivatives involved in \eqref{ns} are regular distributions
and equations \eqref{ns} hold almost everywhere in 
$\OM\times(0,T)$, then $(\n,u)$ is called a strong solution.
\end{definition}

The first main result concerning the global existence and exponential decay of weak solutions is stated as follows.
\begin{theorem}\la{th0}
Assume that the initial data $(\n_0,u_0)$ satisfy
\be\la{wsol1}\ba
0\le \n_0 \in L^\infty,\quad u_0 \in \tilde{H}^1, \quad m_0=\n_0 u_0.
\ea\ee
Then there exists a positive constant $\nu_1$ depending only on 
$\ga$, $\mu$, $E_0$, $\|\n_0\|_{L^\infty}$, $A$, $\OM$ and $\| \na u_0\|_{L^2}$,
such that if $\nu \ge \nu_1 $, the problem \eqref{ns}--\eqref{i3} has at least one weak solution $(\n,u)$ in $\OM \times (0,\infty)$ satisfying 
\be\la{wsol2}\ba
0\le \n(x,t) \leq 2\|\n_0\|_{L^\infty} e^{ \frac{\ga-1}{\ga |\om|} E_0 },
\quad \textnormal{ for any } (x,t) \in \OM \times[0,\infty),
\ea\ee
and for any $0<T<\infty$ and $1 \le p < \infty$,
\be\la{wsol3}\ba
\begin{cases}
\rho\in L^{\infty}(\OM \times (0,\infty)) \cap C([0,\infty);L^p), \\ 
u \in L^2(0,\infty;H^1) \cap L^\infty(0,\infty;H^1), \\
\sqrt{t} u_t \in L^2(0,T;L^2), \sqrt{t} \na u \in L^\infty(0,T;L^p).
\end{cases}
\ea\ee
Moreover, for any $s\in [1,\infty)$, there exist positive constants $C$, $K_0$, and $\widetilde{\nu_0}$ depending on
$s$, $A$, $\OM$, $\ga$,\ $\mu$, $E_0$, $\|\n_0\|_{L^1\cap L^\infty}$, and $\| u_0\|_{H^1}$ such that if $\nu \ge \widetilde{\nu_0}$, then for any $t \ge 1$,
\be\la{wsol4}\ba
\| \n(\cdot,t) - \ol{\n_0} \|_{L^s} + \| \na u(\cdot,t) \|_{L^s} + \| \sqrt{\n} \dot u(\cdot,t) \|_{L^2}
\le C e^{- \alpha_0 t},
\ea\ee
where $\alpha_0=\frac{K_0}{\nu}$.
\end{theorem}

Based on the $\nu$-uniform a priori estimates for the compressible Navier-Stokes equations (\ref{ns}) established in the proof of Theorem \ref{th0}, we can obtain the following incompressible limit result.
\begin{theorem}\la{th01}
Fix the initial data $(\n_0,u_0)$ in $L^\infty\times \tilde{H}^1$ satisfying $\n_0\geq 0$.
Let $\nu_1$ be as in Theorem \ref{th0}.
For each $\nu \ge \nu_1$, let $(\n^{\nu},u^{\nu})$ denote the global weak solution to the problem \eqref{ns}--\eqref{i3} obtained in Theorem \ref{th0}.
Then, as $\nu \to \infty$, the family $(\n^{\nu},u^{\nu})$ has a subsequence converging to a global solution of the following 
inhomogeneous incompressible Navier-Stokes equations:
\be\la{isol1}\ba
\begin{cases}
\n_t+\div(\n u)=0,\\
(\n u)_t+\div(\n u\otimes u) -\mu \Delta u + \na \pi =0, \\
\div u=0, \\
u \cdot n = 0, \quad \curl u=-Au \cdot n^\bot \quad \text{on}\,\,\, \partial\Omega, \\
\end{cases} 
\ea\ee
with initial data $\n(\cdot,0)=\n_0,\ \n u(\cdot,0)=m_0 \triangleq \n_0 u_0$,
and $(\n,u)$ satisfies, for any $1\le p <\infty$ and $0<\tau<T<\infty$,
\be\la{isol2}\ba
\begin{cases}
\rho\in L^{\infty}(\OM \times (0,\infty)) \cap C([0,\infty);L^p), \quad u \in L^2(\OM \times (0,T)), \\
\sqrt{\n} \du, \na \pi, \na^2 u, \na \du \in L^2(\OM \times (\tau,T)), \\
\na u, \sqrt{\n} \du, \na \pi, \na^2 u \in L^\infty(\tau,T;L^2).
\end{cases}
\ea\ee
Moreover, for any $0 < \tau <\infty$, we have
\be\la{isol3}\ba
\div u^{\nu} = O(\nu^{-1/2}) \ \textnormal{ in } L^2(\OM \times (0,\infty)) \cap L^\infty(\tau,\infty;L^2).
\ea\ee
If the initial data $(\n_0,u_0)$ further satisfy
\be\la{insc1}\ba
\div u_0=0,
\ea\ee
then, as $\nu \to \infty$, the whole family $(\n^\nu,u^\nu)$ converges to the unique global solution of \eqref{isol1}.
The limit $(\n,u)$ satisfies, for any $1 \le p < \infty$ and $0<T<\infty$,
\be\la{insc3}\ba
\begin{cases}
\rho\in L^{\infty}(\OM \times (0,\infty)) \cap C([0,\infty);L^p),
\quad \sqrt{\n} u \in C([0,\infty);L^2), \\ 
u\in L^\infty(0,T;H^1), \quad \na u \in L^\infty(0,\infty;L^2), \\
\na u,\ \sqrt{\n} u_t,\ \na^2 u,\ \na \pi \in L^2(\OM \times (0,\infty)), \\ 
\sqrt{t} \na \pi,\ \sqrt{t} \na^2 u \in L^\infty(0,\infty;L^2) \cap L^2(0,T; L^p), \\
\sqrt{\n} u_t \in L^\infty(0,T;L^2) \cap L^2(0,T; L^p), \quad 
\sqrt{t} \na u_t \in L^2(\OM \times (0,T)).
\end{cases}
\ea\ee
\end{theorem}

The following result addresses the global existence and uniqueness of strong solutions.

\begin{theorem}\la{th1}
Suppose that the initial data $(\n_0,u_0)$ satisfy for some $q>2$,
\be \la{ssol1}\ba
0\le \n_0 \in W^{1,q},\quad u_0 \in \tilde{H}^1.
\ea\ee 
Then, for $\nu_1$ as in Theorem \ref{th0}, if $\nu \ge \nu_1 $, the problem \eqref{ns}--\eqref{i3} has a unique strong solution $(\n,u)$ in $\OM \times (0,\infty)$ 
satisfying \eqref{wsol2}
and 
\be\la{ssol4}\ba
\begin{cases}
\rho\in C([0,T];W^{1,q} ), \quad \n_t\in L^\infty(0,T;L^2), \\ 
u\in L^\infty(0,T; H^1) \cap L^{(q+1)/q}(0,T; W^{2,q}), \\ 
t^{1/2}u \in L^2(0,T; W^{2,q}) \cap L^\infty(0,T;H^2), \\
t^{1/2}u_t \in L^2(0,T;H^1), \\
\n u\in C([0,T];L^2), \quad \sqrt{\n} u_t\in L^2(\OM \times(0,T)),
\end{cases} 
\ea\ee
for any $0<T<\infty$.
Moreover, the strong solution $(\n,u)$ satisfies 
\eqref{wsol4}.
\end{theorem}
Next, under higher regularity and compatibility conditions on the initial data $(\n_0,u_0)$,
we can establish the global existence of classical solutions to (\ref{ns}).
\begin{theorem}\la{th2}
Assume that the initial data $(\n_0,u_0)$ satisfy
\be\la{csol1}\ba
0\le \n_0 \in W^{2,q},\quad u_0 \in H^2 \cap \tilde{H}^1,
\ea\ee 
for some $q>2$,
and the following compatibility condition:
\be\la{csol2}\ba
- \mu \Delta u_0 - (\mu + \lm )\nabla \div u_0 +  \nabla P(\n_0)=\n_0^{1/2}g,
\ea\ee
for some $g\in L^2$.
Then, for $\nu_1$ as in Theorem \ref{th0}, if $\nu \ge \nu_1 $, the problem \eqref{ns}--\eqref{i3} has a unique classical solution 
$(\n,u)$ in $\OM \times (0,\infty)$ satisfying \eqref{wsol2}
and 
\be\la{csol4}\ba
\begin{cases}
(\rho,P(\n))\in C([0,T];W^{2,q} ),\quad  (\n_t,P_t)\in L^\infty(0,T;H^1), \\ 
(\n_{tt},P_{tt})\in L^2(0,T;L^2), \\
u\in L^\infty(0,T; H^2) \cap L^2(0,T;H^3), \quad u_t \in L^2(0,T; H^1) \\ 
\na u_t, \na^3u \in L^{(q+1)/q}(0,T;L^q), \\
t^{1/2}\na^3u \in L^\infty(0,T;L^2)\cap L^2(0,T; L^q), \\
t^{1/2}u_t\in L^\infty(0,T;H^1) \cap L^2(0,T;H^2), \\ 
t^{1/2}\na^2(\n u)\in L^\infty(0,T;L^q),\quad \n^{1/2}u_t\in L^{\infty}(0,T;L^2), \\
t\n^{1/2}u_{tt},\quad t\na^2 u_t \in L^{\infty}(0,T;L^2), \\
t\na^3 u \in L^{\infty}(0,T;L^q),\quad t\na u_{tt} \in L^2(0,T;L^2),
\end{cases} 
\ea\ee
for any $0<T<\infty$.
Furthermore, the classical solution $(\n,u)$ satisfies \eqref{wsol4}.
\end{theorem}

Finally, similar to \cite{CL,LX}, we can deduce from $(\ref{wsol4})$
the following large-time behavior of the spatial gradient of the density for the strong solution in Theorem $\ref{th1}$ when vacuum states appear initially.   
\begin{theorem}\la{th3}
In addition to the assumptions in Theorem \ref{th1}, 
we further assume that there exists some point $x_0\in \OM$ such that $\n_0(x_0)=0$. 
Then for any $r>2$, there exists a positive constant $C$ depending only on $r,\ \mu,\ E_0,\ \ga,\ \| \rho_0 \|_{L^1\cap L^\infty}$, such that for any $t \ge 1$
\be\la{pbu0}\ba
\| \na \n(\cdot ,t) \|_{L^r} \ge C e^{\alpha_0 \frac{r-2}{r} t }.
\ea\ee
\end{theorem}

\begin{remark}\la{lrk1}
We conclude from $q>2$ and \eqref{csol4} that
\be\la{csol5}\ba
\n, P(\n) \in C([0,T];W^{2,q})\hookrightarrow C\left([0,T];C^1(\overline{\OM}) \right).
\ea\ee
Furthermore, for any $0<\tau<T<\infty$, the standard embedding implies
\be\la{csol6}\ba
u\in L^\infty(\tau ,T;W^{3,q})\cap H^1(\tau ,T;H^2)\hookrightarrow
C\left([\tau ,T];C^2(\overline{\OM}) \right),
\ea\ee
and	
\be\ba\nonumber
u_t\in L^\infty(\tau ,T;H^2)\cap H^1(\tau ,T;H^1)\hookrightarrow
C\left([\tau ,T]; C(\overline{\OM}) \right).
\ea\ee
By virtue of $(\ref{ns})_1$, \eqref{csol5} and \eqref{csol6}, we have
\be\ba\nonumber
\n_t = -\n \div u - u \cdot \na \n \in C(\overline{\OM} \times [\tau,T]).
\ea\ee
Hence the solution in Theorem \ref{th2} 
is in fact a classical solution 
to the problem \eqref{ns}--\eqref{i3} in $\OM \times (0,\infty)$.
\end{remark}

\begin{remark}\la{lrk2}
We note that Danchin-Mucha \cite[Theorem 2.1]{DM} established the global existence of weak solutions to (\ref{ns}) under the assumptions that the bulk viscosity coefficient is sufficiently large and that $\| \div u_0 \|_{L^2} \le K \nu^{-1/2}$ for some constant $K>0$.
In contrast, our Theorem \ref{th0} demonstrates the global existence of solutions under the sole condition that the bulk viscosity coefficient is sufficiently large, without any extra restrictions on $\div u_0$.
\end{remark}

\begin{remark}\la{lrk3}
Compared with the periodic case $\T^2$ studied by Danchin-Mucha \cite{DM}, to the best of our knowledge, Theorem \ref{th01} appears to be the first global-in-time result on the convergence from compressible Navier-Stokes equations to inhomogeneous incompressible Navier-Stokes equations in bounded domains with Navier-slip boundary conditions.
\end{remark}

\begin{remark}\la{ctoi1}
In contrast to \cite{DM}, we prove that weak solutions of the compressible Navier-Stokes equations \eqref{ns} converge to solutions of the inhomogeneous incompressible Navier-Stokes equations \eqref{isol1} without imposing the condition $\div u_0=0$.
This is consistent with the existence theory in \cite[Theorem 2.1]{L1}, which shows that the inhomogeneous incompressible Navier-Stokes system \eqref{isol1} admits global weak solutions even when $ \div u_0 \neq 0$.
\end{remark}

We now make some comments on the analysis of this paper. 
Note that for initial data satisfying (\ref{csol1})--(\ref{csol2}) with away from vacuum, the local existence and uniqueness of classical solutions to the problem (\ref{ns})--(\ref{i3}) can be established by arguments similar to those in \cite{LLL}.
Therefore, to extend the local classical solution globally in time and allow initial vacuum, it is essential to derive global a priori estimates for smooth solutions to (\ref{ns})--(\ref{i3}) that are independent of the lower bound of the initial density.
In view of the blow-up criterion established in \cite{HLX1}, the key issue is to derive an upper bound for the density that is independent of the lower bound of the initial density.

First, following \cite{DM,HL2,L2}, we introduce the effective viscous flux $G$ (see (\ref{gw}) for its definition) and rewrite $(\ref{ns})_1$ in terms of $G$.
In contrast to the periodic case, the classical commutator theory is not directly applicable in bounded domains, and hence the standard representation of $G$ in terms of Riesz transforms is no longer available.
Nevertheless, we observe that $G$ satisfies an elliptic equation subject to Neumann boundary condition; see (\ref{nsb42}).
As shown in \cite{FLL}, when $\OM$ is a bounded simply connected domain, $G$ can be represented explicitly using the Green function for the unit disk together with a conformal mapping; see (\ref{nsb07}).
Using this representation and arguing as in \cite{DM}, we rewrite $(\ref{ns})_1$ as a transport equation with a linear damping term (see (\ref{nsb91})).
Consequently, a crucial step in deriving the upper bound for the density is to obtain the estimates for $\na u$ and $\sqrt{\n} \dot{u}$.
We first extend the logarithmic interpolation inequality in \cite{DM,D} from the periodic domain to general bounded Lipschitz domains (see Lemma \ref{logbd}).
Combining this inequality with an adaptation of the strategy in \cite{DM} to estimate the convective term (see (\ref{nsb67})), we show that $\log(e+\| \na u \|^2_{L^2})$ can be controlled by $\log(2+\nu)$ (see (\ref{nsb06})).
These key estimates ultimately enable us to derive an upper bound for the density when $\nu$ is sufficiently large (see Proposition \ref{pl1}).

In addition, unlike in the periodic case, boundary terms must be handled carefully in our estimates.
By the Sobolev trace theorem, controlling $\na u$ on the boundary naturally requires control of $\na^2 u$, which is not available at the level of the lower-order a priori estimates.
To overcome this difficulty, we employ the method developed in \cite{CL}.
In particular, the boundary condition $u \cdot n = 0$ on $\p \OM$ yields the identities
\be\la{bjds}\ba
u=(u \cdot n^{\bot})n^{\bot}, \quad (u \cdot \na) u \cdot n=-(u \cdot \na) n \cdot u.
\ea\ee
These identities play a crucial role in estimating the boundary terms.
Once the upper bound for the density has been established, arguments similar to those in \cite{CL,FLL,LLL,HLX3} yield the exponential decay and higher-order derivative estimates for the solution, thereby allowing the local solution to be extended globally in time.

Finally, we study the singular limit as $\nu \to \infty$, in which solutions of the compressible Navier-Stokes equations converge to solutions of the inhomogeneous incompressible Navier-Stokes equations.
Following the approach in \cite{DM}, a crucial step in establishing this convergence is to derive a $\nu$-uniform bound for $\sqrt{\n} \dot{u}$.
Using the logarithmic interpolation inequality together with standard energy estimates, we obtain, for any $0<\tau<T<\infty$, a $\nu$-uniform bound for $\sqrt{\n} \dot{u}$ in $L^2(\OM \times (\tau,T))$.
Based on this estimate, we further derive, for any $0<\tau<T<\infty$, $\nu$-uniform bounds for $\sqrt{\n} \dot{u}$ in $L^\infty(\tau,T;L^2(\OM))$ and for $\na \dot{u}$ in $L^2(\OM \times (\tau,T))$.
These $\nu$-uniform estimates, together with standard compactness arguments, imply that the family of solutions to the compressible Navier-Stokes has a subsequence converging to a global solution of the inhomogeneous incompressible Navier-Stokes equations.
Moreover, motivated by \cite{PT}, we prove the uniqueness of weak solutions to the inhomogeneous incompressible Navier-Stokes equations subject to the Navier-slip boundary condition.
Combining this uniqueness result with the $\nu$-uniform estimates yields the convergence of the whole family in the incompressible limit.

The rest of this paper is organized as follows: Section 2 introduces some essential inequalities and known facts.
In Section 3, we establish a time-uniform upper bound for the density that is independent of the lower bound of the initial density.
Building on this density estimate, Section 4 is devoted to deriving the necessary higher-order derivative estimates.
Finally, Section 5 presents the proofs of our main results, Theorems \ref{th0}--\ref{th3}.

\section{Preliminaries}
In this section, we will recall some known facts and elementary inequalities
which will be used frequently later.

First, we have the following local existence theory of the classical solution, 
which can be proved in a manner similar to \cite{LLL}.
\begin{lemma}\la{lct}
Assume $\left( \n _0 ,u_0  \right)$ satisfies that for some $q>2$
\be\ba\nonumber
\n _0 \in W^{2,q}, \quad \inf\limits_{x\in\OM}\n_0(x) >0, \quad u_0 \in H^2 \cap \tilde{H}^1,
\ea\ee
and the compatibility condition \eqref{csol2}. 
Then there is a small time $T>0$ and a constant $C_0>0$ both depending only on
$\mu,\ \lambda,\ \ga,\ q,\ \| \n_0\|_{W^{2,q}},\ \|u_0\|_{H^2},\  \inf\limits_{x\in\OM}\n_0(x)$ and $\|g\|_{L^2}$,  
 such that there exists a unique classical solution $(\n,u)$ to the problem 
 $(\ref{ns})-(\ref{i3})$ in $\OM \times (0,T]$ satisfying (\ref{csol4})
and
\be\ba\nonumber
\inf\limits_{(x,t)\in \OM \times (0,T)}\n(x,t) \ge C_0 >0.
\ea\ee

\end{lemma}
Next, the following Gagliardo-Nirenberg's inequalities (see \cite{TG}) will be used frequently later.
\begin{lemma}\la{gn1}
Let $u\in H^1(\OM)$. For any $2<p<\infty$, there exists a positive constant C depending only on $\OM$ such that 
\be\ba\la{gn11}
\| u\|_{L^p} \le Cp^{1/2}\| u\|^{2/p}_{L^2} \|  u\|^{1-2/p}_{H^1}.
\ea\ee
Furthermore, $\| u\|_{H^1}$ can be replaced by $\| \na u\|_{L^2}$ provided
\bnn
u \cdot n |_{\p \OM}=0 \  or \  \int_\OM u dx =0.
\enn

Moreover, for $1\le r <\infty$, $2<q<\infty$, there exists a positive constant $C$ depending only on $r$, $q$, and $\OM$ such that for every function $v\in W^{1,q}(\OM)$, it holds that
\be\ba\la{gn12}
\|v-\ol{v}\|_{L^\infty} \le C\|v-\ol{v}\|^{r(q-2)/2q+r(q-2)}_{L^r} \|\na v\|^{2q/2q+r(q-2)}_{L^q}.
\ea\ee
\end{lemma}

The following Poincar\'e-type inequality can be found in \cite{F}.
\begin{lemma}\la{pt}
Let $\OM \subset \rr$ be a bounded Lipschitz domain.
Assume that $v\in H^1(\OM)$ and that $\n$ is a non-negative function satisfying
\be\ba\nonumber
0 < M_1\leq \int \n dx, \quad \int \n^r dx \leq M_2,
\ea\ee
for some $r>1$.
Then there exists a positive constant $C$ depending only on $M_1,\ M_2$ and, $r$ such that
\be\ba\la{pt1}
\|v\|_{L^2}^2 \leq C\int \n |v|^2 dx + C \|\na v\|_{L^2}^2.
\ea\ee
\end{lemma}

The following div-curl estimates will be used frequently; see \cite{AJ,MD,WWV}.
\begin{lemma}\la{dc}
Let $k \ge 0$ be an integer and let $1<p<\infty$.
Assume that $\OM$ is a simply connected bounded domain in $\rr$ with $C^{k+1,1}$ boundary $\p \OM$.
Then there exists a positive constant $C$ depending only on $k$, $p$, and $\OM$ such that for every $u\in W^{k+1,p}$ with $u \cdot n=0$ on $\p \OM$,
\be\ba\la{dc1}
\| u\|_{W^{k+1,p}} \le C\left( \|\div u\|_{W^{k,p}} +\| \curl u \|_{W^{k,p}} \right).
\ea\ee
\end{lemma}

More generally, we have the following weighted div-curl estimates; see \cite{FLL,FLW}.
\begin{lemma}\la{wdc}
Let $\OM$ be a simply connected bounded domain in $\rr$ with Lipschitz boundary $\p \OM$.
Then there exist positive constants $C$ and $\hat{\de}$ depending only on $\OM$ such that for any $\de \in (0,\hat{\de})$ and any $u\in H^2(\OM)$ with $u \cdot n=0$ on $\p \OM$,
\be\ba\la{wdc1}
\int_\OM |u|^{\de} |\na u|^2 dx \le C \int_\OM |u|^{\de} \left( (\div u)^2 +(\curl u)^2  \right)dx.
\ea\ee
\end{lemma}

Subsequently, the following estimates on the material derivative of $u$ play an important role in the higher order estimates, 
whose proof can be found in \cite[Lemma 4.1]{FLL}.
\begin{lemma}\la{emdu}
For any $1\le p<\infty$, there exist two positive constants $\Lambda_1$ and $\Lambda _2$,
where $\Lambda _1$ depends on $p$ and $\OM$, while $\Lambda _2$ depends only on $\OM$,
such that
\be\la{emdu1}\ba
\| \dot{u} \|_{L^p} \le \Lambda_1 \left( \| \na \dot{u} \|_{L^2}+\| \na u \|^2_{L^2} \right),
\ea\ee
\be\la{emdu2}\ba
\| \na \dot{u} \|_{L^2} \le \Lambda_2 \left( \| \div \dot{u} \|_{L^2}+\| \curl \dot{u} \|_{L^2}+\| \na u \|^2_{L^4} \right).
\ea\ee
\end{lemma}

To estimate $\| \na u\|_{L^{\infty}}$ and $\| \na \n\|_{L^{q}}$ 
we require the following Beale-Kato-Majda type inequality, 
which was established in \cite{K} when $\div u \equiv 0$. 
For further reference, we direct readers to \cite{BKM,CL,HLX1}.
\begin{lemma}\la{bkm}
Let $2<q<\infty$.
Assume that $\OM$ is a simply connected bounded $C^{3,1}$ domain in $\rr$ and $u \in \left\{ W^{2, q}(\Omega ) \big| u \cdot n=0, \curl u = -A u \cdot n^\bot \,\,\,\text{on}\,\,\, \partial\Omega \right\}$.
Then there exists a positive constant $C$ depending only on $\OM$ and $q$ such that
\be\ba\nonumber
 \|\na u\|_{L^\infty} \le C \left( \|\div u\|_{L^\infty}+ \|\o\|_{L^\infty} \right) \log \left(e+ \|\na^2 u\|_{L^q} \right)+ C\|\na u\|_{L^2}+C.
\ea\ee
\end{lemma}

The following logarithmic interpolation inequality in $\rr$ can be found in \cite[Lemma 2.4]{Z}.
\begin{lemma}\la{logrr}
Suppose that $\n \in L^\infty(\rr)$ satisfies $0 \le \n \le \rs $ and that $u \in H^1(\rr)$.
Then there exists a positive constant $C$ depending only on $\rs$ such that
\be\la{logrr1}\ba
\| \sqrt{\n} u \|^2_{L^4(\rr)} \le C \left( 1+\| \sqrt{\n} u \|_{L^2(\rr)} \right)
\| u \|_{H^1(\rr)} \log^{\frac{1}{2}} \left( 2+\| u \|^2_{H^1(\rr)} \right).
\ea\ee
\end{lemma}

Next, we extend the preceding inequality to bounded domains.
\begin{lemma}\la{logbd}
Let $\OM \subset \rr$ be a bounded Lipschitz domain.
Assume that $\n \in L^\infty(\OM)$ satisfies $0 \le \n \le \rs $ and that $u \in H^1(\OM)$.
Then there exists a positive constant $C$ depending only on $\rs$ and $\OM$ such that
\be\la{logbd1}\ba
\| \sqrt{\n} u \|^2_{L^4(\OM)} \le C \left( 1+\| \sqrt{\n} u \|_{L^2(\OM)} \right)
\| u \|_{H^1(\OM)} \log^{\frac{1}{2}} \left( 2+\| u \|^2_{H^1(\OM)} \right).
\ea\ee
\end{lemma}
\begin{proof}
First, we extend $\n$ by zero outside $\OM$ and denote this extension by $\hat{\n}$.
Then $\hat{\n} \in L^\infty(\rr)$ and $0 \le \hat{\n} \le \rs $.

Fix a bounded open set $V$ such that $\ol{\OM} \subset V$.
By the Sobolev extension theorem \cite[Chapter 5]{EL}, for $u \in H^1(\OM)$ there exists a function $\hat{u} \in H^1(\rr) $ such that
\be\la{logbd2}\ba
\hat{u}=u \ \ \text{ a.e. in } \  \OM,
\quad \operatorname{supp} \hat{u} \subset V,
\quad \| \hat{u} \|_{H^1(\rr)} \le C \| u \|_{H^1(\OM)}, 
\ea\ee
where $C$ depends only on $\OM$ and $V$.

Applying Lemma \ref{logrr} to $(\hat{\n},\hat{u})$, we obtain
\be\la{logbd3}\ba
\| \sqrt{\hat{\n}} \hat{u} \|^2_{L^4(\rr)}
\le C \left( 1+\| \sqrt{\hat{\n}} \hat{u} \|_{L^2(\rr)} \right)
\| \hat{u} \|_{H^1(\rr)} \log^{\frac{1}{2}} \left( 2+\| \hat{u} \|^2_{H^1(\rr)} \right).
\ea\ee
Combining this with (\ref{logbd2}) and the definition of $\hat{\n}$ yields (\ref{logbd1}), which completes the proof.
\end{proof}

Finally, we recall the following ``inversion'' operator of the divergence on bounded domains (see \cite{CL} for details).
\begin{lemma}\la{iod}
Let $\OM \subset \rr$ be a bounded Lipschitz domain.
For $1<p<\infty$, there exists a bounded linear operator
\be\ba\nonumber
\mathcal{B} : \left\{ f\in L^p(\OM) : \ \int_\Omega f dx=0 \right\}&\rightarrow (W^{1,p}_0(\OM))^2,
\ea\ee
such that $v=\mathcal{B}(f)$ satisfies
\be\la{iod1}\ba
\begin{cases}
\mathrm{div}v=f&\ \textnormal{ in }\Omega,\\
v=0&\ \textnormal{ on }\partial\Omega.
\end{cases}
\ea\ee
Moreover, the following estimates hold.

\textup{(1)} For $1<p<\infty$, there is a positive constant $C$ depending on $\Omega$ and $p$, such that
\bnn
\|\mathcal{B}(f)\|_{W^{1,p}}\leq C(p)\|f\|_{L^p}.
\enn

\textup{(2)} If $f=\mathrm{div}h$ for some $h\in L^p$ with $h\cdot n=0$ on $\partial\Omega$, then $v=\mathcal{B}(f)$ is a weak solution of \eqref{iod1} and
\bnn
\|\mathcal{B}(f)\|_{L^p}\leq C(p) \|h\|_{L^p}.
\enn
\end{lemma}

\section{A Priori Estimates (I): Upper Bound of the density}
In this section, we assume that the initial data $(\n_0,u_0)$ satisfy (\ref{wsol1}), (\ref{ssol1}), (\ref{csol1}), and (\ref{csol2}) with $\n_0>0$.
Let $T>0$ be fixed, and let $(\n,u)$ be a local classical solution to the problem (\ref{ns})--(\ref{i3}) on  $\OM \times (0,T]$, whose existence is guaranteed by Lemma \ref{lct}.

We introduce the effective viscous flux $G$ defined by:
\be\ba\la{gw}
G \triangleq (2\mu + \lam)\div u - (P-\ol{P}).
\ea\ee

We also set
\be\ba\nonumber
B^2_1(t) \triangleq \int \left( \nu (\div u)^2 + \mu \o^2 + \frac{1}{\nu} (P-\ol{P})^2 \right) dx
+ \mu \int_{\p \OM} A|u|^2 ds.
\ea\ee
and
\be\ba\nonumber
B_2^2(t)\triangleq\int \rho|\dot{u}|^2dx.
\ea\ee

This section aims to derive the following a priori estimates.
\begin{proposition}\la{pl1}
There exists a positive constant $\nu_1$ depending only on
$\mu$, $\gamma$, $\OM$, $A$, $E_0$, $\| \n_0 \|_{L^\infty}$, and $\| u_0 \|_{H^1}$
such that if $(\n,u)$ satisfies
\be\ba\la{p11}
\sup_{0\leq t\leq T} \|\n\|_{L^\infty} 
\leq 2 \| \n_0 \|_{L^\infty} e^{ \frac{\ga-1}{\ga |\OM|} E_0 },
\ea\ee
then
\be\ba\la{p12}
\sup_{0\leq t\leq T} \|\n\|_{L^\infty} 
\leq \frac{3}{2} \| \n_0 \|_{L^\infty} e^{ \frac{\ga-1}{\ga |\OM|} E_0 },
\ea\ee
provided $\nu \geq \nu_1$.
\end{proposition}

The proof of Proposition \ref{pl1} will be postponed to the end of this subsection.

We first state the standard energy estimate.
\begin{lemma}\la{l1}
Suppose that $(\n,u)$ is a classical solution to \eqref{ns}--\eqref{i3} on $\OM \times (0,T]$.
Then we have:
\be\ba\la{nsb1}
& \sup_{0\leq t\leq T} \int \left( \frac{1}{2}\rho |u|^2 + \frac{P}{\ga-1} \right) dx
+ \int_0^T \left( \nu \| \div u \|^2_{L^2} + \mu \| \o \|^2_{L^2} \right) dt \\
& \quad + \mu \int_0^T \int_ {\p\OM} A |u|^2 dsdt
\le E_0,
\ea\ee
where $E_0$ is defined by \eqref{e0} and 
$\o$ is defined by \eqref{xd}.
\end{lemma}
\begin{proof}
Multiplying $(\ref{ns})_2$ by $u$ and integrating the resulting equation over $\OM$, after using the boundary condition (\ref{i3}) and $(\ref{ns})_1$, we obtain (\ref{nsb1}).
\end{proof}

\begin{lemma}\la{l2}
Let $(\n,u)$ be a classical solution to \eqref{ns}--\eqref{i3} satisfying \eqref{p11}.
Then there exists a positive constant $C$ depending only on
$\mu$, $\gamma$, $\| \n_0 \|_{L^\infty}$, $E_0$, and $\OM$ such that
\be\la{nsb02}\ba
\int_0^T \int(P-\ol{P})^{2} dxdt\le C \nu.
\ea\ee
\end{lemma}
\begin{proof}
Since $P$ satisfies
\be\la{nsb22}\ba
P_t+\div(Pu)+(\gamma-1)P\div u=0,
\ea \ee
integrating this equation over $\OM$ and using (\ref{i3}), we arrive at
\be\la{nsb23}\ba
\ol{P_t}+(\ga-1)\ol{P\div u}=0.
\ea\ee
Multiplying $(\ref{ns})_2$ by $\mathcal{B}[P-\ol{P} ]$, integrating over $\Omega$,
and applying (\ref{pl1}), (\ref{nsb1}), and H\"older's inequality, we conclude that
\be\la{nsb24} \ba
\int(P-\ol{P} )^2 dx 
&= \left(\int\rho u\cdot \mathcal{B}[P-\ol{P}] dx\right)_t- \int\rho u\cdot\mathcal{B}[P_t-\ol{P_t}]dx  \\
& \quad -\int\rho u\cdot\nabla\mathcal{B}[P-\ol{P}]\cdot udx +\mu \int \p_i u \cdot \p_i \mathcal{B}[P-\ol{P}] dx \\
& \quad  + (\mu+\lambda)\int(P-\ol{P})\div udx \\
& \le \left(\int\rho u\cdot\mathcal{B}[P-\ol{P}] dx\right)_t+\|\n  u\|_{L^2}\|\mathcal{B}[P_t-\ol{P_t} ]\|_{L^2} \\
& \quad +C \| \n \|_{L^4} \| u\|_{L^4}^{2} \|P-\ol P\|_{L^4} +C\|P-\ol P\|_{L^2} \left( \|\nabla u\|_{L^2}+ \nu \| \div u \|_{L^2}  \right) \\
& \leq \left(\int\rho u\cdot\mathcal{B}[P-\ol{P}] dx\right)_t+\frac{1}{2} \|P-\ol{P}\|_{L^2}^2 +C \left( \|\nabla u\|^2_{L^2}+ \nu^2 \| \div u \|^2_{L^2}  \right),
\ea\ee 
where in the last inequality we have used the following simple fact:
\bnn\ba \|\mathcal{B}[P_t-\ol{P_t} ]\|_{L^2}
&=\|\mathcal{B} [\div(Pu)] + (\ga-1) \mathcal{B} [P\div u-\ol{P\div u}]\|_{L^2} \\
& \le C\left( \| P u\|_{L^2}+ \| P \div u\|_{L^2} \right) \\
&\le C \|\na u\|_{L^2},
\ea\enn
due to (\ref{nsb22}), (\ref{nsb23}) and Lemma \ref{iod}.

Integrating (\ref{nsb24}) over $(0,T)$ and using (\ref{pl1}), (\ref{nsb1}), and Lemma \ref{iod},
we derive (\ref{nsb02}) and finish the proof of Lemma \ref{l2}.
\end{proof}

\begin{lemma}\la{l3}
Let $(\n,u)$ be a classical solution to \eqref{ns}--\eqref{i3} satisfying \eqref{p11}.
Then there exists a positive constant $C$ depending only on
$\mu$, $\gamma$, $\| \n_0 \|_{L^\infty}$, $E_0$, and $\OM$ such that
\be\la{nsb03}\ba
\sup_{0\le t\le T}\int \n |u|^{2+\de}dx\le C \nu,
\ea\ee
where
\be\la{nsb003}\ba
\de \triangleq \nu^{-\frac{1}{2}} \de_0,
\ea\ee
for some positive constant $\de_0 \le \frac{1}{2} \mu^{\frac{1}{2}}$ depending only on $\mu$ and $\OM$.
\end{lemma}
\begin{proof}
First, we multiply $(\ref{ns})_2$ by $|u|^\de u$ and integrate the resulting equation over $ \OM$ to obtain
\be\la{nsb31}\ba 
& \frac{1}{(2+\de)} \frac{d}{dt}\int \n |u|^{2+\de}dx + \int|u|^\de \left(\mu |\o|^2+ \nu(\div u)^2\right) dx +\mu \int_{\p \OM} A|u|^{2+\de}dS \\
& \le C \de  \int  \left(\nu |\div u|+\mu |\o| \right)  |u|^\de |\na u|dx +C \int |P-\ol{P}| |u|^\de |\na u|dx \\
& \triangleq I_1+I_2.
\ea\ee
For $I_1$, it follows from (\ref{wdc1}) and Cauchy's inequality that
\be\la{nsb32}\ba
I_1 &\le \frac{1}{2}\int|u|^\de \left(\mu |\o|^2+ \nu(\div u)^2\right) dx +\frac{C\de^2 \nu}{2} \int |u|^\de |\na u|^2 dx \\
& \le \frac{1 + C_1 \de^2 \nu}{2}\int|u|^\de \left(\mu |\o|^2+ \nu(\div u)^2\right) dx,
\ea\ee
provided $\de \in (0,\hat{\de})$, where $C_1$ depends only on $\mu$ and $\OM$.

For $I_2$, making use of Young's and Poincar\'e's inequalities, we obtain
\be\la{nsb33}\ba
I_2 &\le C \int |P-\ol{P}| \left( 1 + |u| \right)  |\na u| dx \\
& \le C \| P-\ol{P} \|_{L^2} \| \na u\|_{L^2} + C \| u \|_{L^2} \| \na u\|_{L^2} \\
& \le C \left( \| P-\ol{P} \|^2_{L^2} + \| \na u\|^2_{L^2} \right).
\ea\ee
Putting (\ref{nsb32}) and (\ref{nsb33}) into (\ref{nsb31}) and taking 
$\de_0= \min \left\{ \frac{1}{2} \sqrt{\mu},\sqrt{\mu}\hat{\de},\frac{1}{\sqrt{2 C_1}} \right\}$, we get
\be\la{nsb35}\ba 
\frac{d}{dt}\int \n |u|^{2+\de}dx \le C \left( \| P-\ol{P} \|^2_{L^2} + \| \na u\|^2_{L^2} \right).
\ea\ee
Moreover, we deduce from Poincar\'e's inequality that
\be\la{nsb36}\ba
\int \n_0 |u_0|^{2+\de} dx \le C \| u_0 \|^{2+\de}_{H^1} 
\le C \| \na u_0 \|^{2+\de}_{L^2} \le C.
\ea\ee
Integrating (\ref{nsb35}) over $(0,T)$ and using (\ref{nsb36}), (\ref{nsb1}), and (\ref{nsb02}), we derive (\ref{nsb03}) and complete the proof of Lemma \ref{l3}.
\end{proof}

For $2<p<\infty$, the following estimate of $\| \na u \|_{L^p}$ will be frequently used and is crucial in the subsequent analysis.
\begin{lemma}\la{l4}
Let $(\n,u)$ be a classical solution to \eqref{ns}--\eqref{i3} satisfying \eqref{p11}.
Then, for any $2 \le p<\infty$, there exists a positive constant $C$ depending only on
$p$, $\mu$, $\gamma$, $\| \n_0 \|_{L^\infty}$, $E_0$, and $\OM$ such that
\be\la{nsb04}\ba
\| \nabla u\|_{L^{p}} \le C B^{\frac{2}{p}}_1 B^{1-\frac{2}{p}}_2
+ CB_1 + \frac{C}{\nu} \| P-\ol{P} \|_{L^p}.
\ea\ee
\end{lemma}
\begin{proof}
First, we rewrite $(\ref{ns})_2$ as 
\be\la{nsb41}\ba
\n\dot{u}= \na G + \mu\na^{\bot}\o,
\ea\ee
which together with the boundary condition (\ref{i3}) implies that $G$ and $\o$ satisfy the following elliptic equations, respectively:
\be\la{nsb42}\ba
\begin{cases}
\Delta G=\div \left(\rho \dot{u} \right) &\mathrm{in}\, \,  \OM, \\
\frac {\p G}{\p n}=\left(\rho \dot{u}-\mu \na^{\bot} \o \right) \cdot n &\mathrm{on}\, \,  \p \OM,
\end{cases}
\ea\ee
and
\be\la{nsb43}\ba
\begin{cases}
\mu\Delta\o =\na ^\bot \cdot \left(\rho \dot{u} \right)& \mbox{ in } \OM, \\
\o=-Au \cdot n^{\bot}& \mbox{ on } \p\OM.  
\end{cases}
\ea\ee

Based on the standard $L^p$ estimate of elliptic equations as stated in \cite[Lemma 4.27]{NS}, we obtain that for any integer $k \ge 0$ and $1<p<\infty$,
\be\la{nsb44}\ba
\| \na G \|_{W^{k,p}}+\| \na \o \|_{W^{k,p}} \le C\left(\| \n \dot u\|_{W^{k,p}}
+ \| \na u\|_{W^{k,p}} \right),
\ea\ee
where $C$ depends only on $p,\ k,\ \mu,\ A,\ \OM$.

In particular, by making use of H\"older's and Poincar\'e's inequalities,
we conclude that
\be\la{nsb45}\ba
\| G \|_{H^1}+\| \o \|_{H^1} 
\le C\left(\| \n \dot u\|_{L^2} + \| \na u\|_{L^2} \right).
\ea\ee
Consequently, we deduce from (\ref{gn11}), (\ref{dc1}), (\ref{pl1}), (\ref{nsb45}), and Young's inequality that
\be\la{nsb46}\ba
\| \na u \|_{L^p} & \le C \left(\| \div u\|_{L^p}+\| \o\|_{L^p} \right) \\
& \le C \left( \frac{1}{\nu} \| G \|_{L^p}+ \frac{1}{\nu} \| P-\ol{P} \|_{L^p} + \| \o\|_{L^p} \right) \\
& \le C \left( \frac{1}{\nu} \| G \|^{\frac{2}{p}}_{L^2} \| G \|^{1-\frac{2}{p}}_{H^1}
+ \frac{1}{\nu} \| P-\ol{P} \|_{L^p} + \| \o \|^{\frac{2}{p}}_{L^2} \| \o \|^{1-\frac{2}{p}}_{H^1} \right) \\
& \le C B^{\frac{2}{p}}_1 \left(\| \n \dot u\|_{L^2} + \| \na u\|_{L^2} \right)^{1-\frac{2}{p}} 
+\frac{C}{\nu} \| P-\ol{P} \|_{L^p} \\
& \le C B^{\frac{2}{p}}_1 B^{1-\frac{2}{p}}_2
+ C B_1 + \frac{C}{\nu} \| P-\ol{P} \|_{L^p},
\ea\ee
which implies (\ref{nsb04}) and completes the proof of Lemma \ref{l4}.
\end{proof}

\begin{lemma}\la{l6}
Let $(\n,u)$ be a classical solution to \eqref{ns}--\eqref{i3} satisfying \eqref{p11}.
Then there exists a positive constant $C$ depending only on
$\mu$, $\gamma$, $A$, $\| \n_0 \|_{L^\infty}$, $E_0$, and $\OM$ such that
\be\ba\la{nsb06}
\sup_{0\leq t\leq T} \log \left( 2 + B^2_1(t) \right) + \int_0^T \frac{B_2^2}{ 2 + B^2_1(t) } dt
\le C \log \left( 2 + B^2_1(0) \right),
\ea\ee
and, for any $\tau \in (0,T)$,
\be\ba\la{nsb006}
\sup_{\tau \le  t \le T} \left( \| \na u \|_{L^2}^2 + \nu \| \div u \|_{L^2}^2 \right)
+ \int_\tau^T \| \sqrt{\n} \dot{u} \|_{L^2}^2 dt \le C + C \tau^{-C}.
\ea\ee
\end{lemma}
\begin{proof}
First, multiplying both sides of $(\ref{ns})_2$ by $u_t$ and integrating
the resulting equality over $\OM$, we derive
\be\ba\nonumber
& \frac{d}{dt} \left( \frac{\nu}{2}\|\div u\|_{L^2}^2 + \frac{\mu}{2}
\|\o\|^2_{L^2} \right)
+\frac{\mu}{2} \frac{d}{dt} \int_{\p \OM} A|u|^2 ds	+\int\n|\dot{u}|^2dx \\
&= \int (P-\ol{P}) \div u_t dx + \int \n\du\cdot(u\cdot\na)udx\\
& \leq \frac{1}{2} \int \n |\du|^2 dx 
+\frac{d}{dt} \int (P-\ol{P}) \div u dx- \int P_t \div u dx 
+ \frac{1}{2}\int \n|u|^2|\na u|^2 dx \\
& \leq \frac{1}{2} \int \n |\du|^2 dx 
+\frac{d}{dt} \int (P-\ol{P}) \div u dx- \frac{1}{\nu}\int P_t G dx
-\frac{1}{\nu}\int P_t (P-\ol{P}) dx \\
& \quad + \frac{1}{2}\int \n|u|^2|\na u|^2 dx \\
& =  \frac{d}{dt} \left(\int (P-\ol{P})\div u dx
-\frac{1}{2\nu}\| P-\ol{P} \|^2_{L^2} \right) 
+ \frac{1}{2} \int \n |\du|^2 dx - \frac{1}{\nu} \int P_t G dx \\
& \quad + \frac{1}{2}\int \n|u|^2|\na u|^2 dx,
\ea\ee
where we have used $(\ref{ns})_1$, (\ref{i3}) and (\ref{gw}).

Therefore, we have
\be\ba\la{nsb62}
\frac{d}{dt}B_3(t) + \frac{1}{2} B_2^2
&
\leq - \frac{1}{\nu} \int P_t G dx +\frac{1}{2} \int \n|u|^2|\na u|^2 dx \\
&
\triangleq I_1 + I_2,
\ea\ee
with
\be\ba\la{nsb63}
B_3(t) \triangleq \frac{\nu}{2}\|\div u\|_{L^2}^2 + \frac{\mu}{2}\|\o\|^2_{L^2} + \frac{\mu}{2} \int_{\p \OM} A|u|^2 ds
+\frac{1}{2\nu}\| P-\ol{P} \|^2_{L^2} 
-\int( P-\ol{P} )\div udx.
\ea\ee

Next we proceed to estimate $I_1$ and $I_2$ sequentially.
We deduce from (\ref{nsb22}), (\ref{gw}), (\ref{nsb45}) and Poincar\'e's inequality that
\be\ba\la{nsb66}
I_1
& = -\frac{1}{\nu}\int Pu\cdot\na G dx 
+ \frac{\ga-1}{\nu} \int P\div u Gdx\\
& \leq \frac{C}{\nu} \|u\|_{L^2} \|\na G\|_{L^2}+\frac{C}{\nu} \|\div u\|_{L^2} \|G\|_{L^2} \\
& \leq \frac{C}{\nu} \| \na u \|_{L^2} \left(\| \n \dot u\|_{L^2}
+ \| \na u\|_{L^2} \right) \\
& \le \frac{1}{8} B_2^2 + C B_1^2.
\ea\ee
Then, we turn to estimating $I_2$. By virtue of (\ref{logbd1}), (\ref{nsb04}) and H\"older's inequality, we arrive at
\be\ba\la{nsb67}
|I_2|
& \leq C\| \sqrt{\n} u \|^2_{L^4} \|\na u\|_{L^4}^2 \\
& \leq C \| \sqrt{\n} u \|^2_{L^4}
\left( B_1 B_2 + B^2_1 + \frac{1}{\nu^2} \| P-\ol{P} \|^2_{L^4} \right) \\
& \leq \frac{1}{8} B^2_2 + C \left( B^2_1 \| \sqrt{\n} u \|^4_{L^4}+\| \sqrt{\n} u \|^4_{L^4}+\frac{1}{\nu^4} \| P-\ol{P} \|^4_{L^4} +B^2_1 \right) \\
& \le \frac{1}{8} B^2_2 + C B^2_1 + C (1+B^2_1)
\left( 1+\| \sqrt{\n} u \|^2_{L^2} \right) \| u \|^2_{H^1} \log \left( 2+\| u \|^2_{H^1} \right) \\
& \le \frac{1}{8} B^2_2 + C B^2_1 + C (1+B^2_1) B^2_1 \log \left( 2+B^2_1 \right).
\ea\ee

Substituting (\ref{nsb66}) and (\ref{nsb67}) into (\ref{nsb62}), we obtain
\be\ba\la{nsb69}
\frac{d}{dt}B_3(t) + \frac{1}{4} B_2^2
& \le C B^2_1 + C (1+B^2_1) B^2_1 \log \left( 2+B^2_1 \right).
\ea\ee

Moreover, by using H\"older's and Young's inequalities, we can derive that there 
exists a positive constant $\check{C}$ depending only on $\ga,\ \mu,\ \rs$ such that
\be\ba\nonumber
\left| \int(P-\ol{P})\div udx \right|
\le \frac{\nu}{4} \|\div u\|_{L^2}^2 + \check{C}.
\ea\ee
We set
\be\ba\nonumber
B_4(t) \triangleq B_3(t) + \check{C},
\ea\ee
which together with (\ref{nsb63}) implies that
\be\ba\la{nsb612}
\frac{1}{4} B^2_1(t) \leq B_4(t) \leq 2 \left( B^2_1(t) + \check{C} \right).
\ea\ee
Consequently, we conclude from (\ref{nsb69}) and (\ref{nsb612}) that
\be\ba\la{nsb613}
\frac{d}{dt} \left( 2 + B_4(t) \right) + \frac{1}{4} B_2^2
& \leq C B^2_1 + C \left( 2 + B_4 \right) B^2_1 \log \left( 2 + B_4 \right).
\ea\ee
Dividing (\ref{nsb613}) by $2 + B_4(t)$, we arrive at
\be\ba\la{nsb614}
\frac{d}{dt} \log \left( 2 + B_4(t) \right) + \frac{B_2^2}{4(2 + B_4(t))} 
& \leq C B^2_1 + C B^2_1 \log \left( 2 + B_4 \right).
\ea\ee
Using (\ref{nsb1}) and (\ref{nsb02}), we have
\be\la{nsb615}\ba
\int_0^T B_1^2 dt \le C \int_0^T \left( \nu \| \div u \|_{L^2}^2 + \| \na u \|_{L^2}^2 + \frac{1}{\nu} \| P - \ol{P} \|_{L^2}^2 \right) dt \le C.
\ea\ee
Applying Gr\"onwall's inequality to (\ref{nsb614}) and using (\ref{nsb615}), we obtain
\be\ba\nonumber
\sup_{0 \le t \le T} \log \left( 2 + B_4(t) \right) + \int_0^T \frac{B_2^2}{4(2 + B_4(t))} dt
\le C \log \left( 2 + B_4(0) \right) + C,
\ea\ee
which together with (\ref{nsb612}) implies (\ref{nsb06}).

On the other hand, for any $\tau \in (0,T)$, by (\ref{nsb615}), there exists $\tau_* \in (\frac{\tau}{2},\tau)$ such that
\be\la{nsb616}\ba
B_1^2(\tau_*) \le \frac{2}{\tau} \int_{\frac{\tau}{2}}^{\tau} B_1^2 dt \le \frac{C}{\tau}.
\ea\ee
Applying Gr\"onwall's inequality to (\ref{nsb614}) over $(\tau_*,T)$ and using (\ref{nsb616}), we get
\be\ba\nonumber
& \sup_{\tau_* \le t \le T} \log \left( 2 + B_4(t) \right) + \int_{\tau_*}^T \frac{B_2^2}{4(2 + B_4(t))} dt \\
& \le C \log \left( 2 + B_4(\tau_*) \right)
\le C \log \left( 2 + 2B^2_1(\tau_*) + 2\check{C} \right) \\
& \le C \log \left( C + C \tau^{-1} \right),
\ea\ee
which together with $\tau_*<\tau$ yields
\be\la{nsb617}\ba
\sup_{\tau \le  t \le T} B_4(t) + \int_\tau^T B_2^2 dt \le  C + C \tau^{-C}.
\ea\ee
Combining (\ref{nsb617}) and (\ref{nsb612}), we arrive at
\be\la{nsb618}\ba
\sup_{\tau \le  t \le T} \left( \| \na u \|_{L^2}^2 + \nu \| \div u \|_{L^2}^2 \right)
+ \int_\tau^T \| \sqrt{\n} \dot{u} \|_{L^2}^2 dt \le C + C \tau^{-C}.
\ea\ee
This gives (\ref{nsb006}) and completes the proof of Lemma \ref{l6}.
\end{proof}

Next, we use the above a priori estimates to obtain the upper bound of $\n$.

First, in light of $(\ref{ns})_1$ and (\ref{gw}), we have
\be\ba\nonumber
\frac{D}{Dt} \log \n+\frac{1}{\nu}(P-\ol{P}) =-\frac{1}{\nu}G.
\ea\ee

Consequently, to derive the upper-bound estimate of $\n$, we need to obtain the $L^\infty$ estimate of $G$.
To achieve this objective, we employ the method from \cite{FLL} to acquire the point-wise representation of $G$.

It follows from (\ref{nsb42}) that for any $t\in [0,T]$, $G$ satisfies the following Neumann problem:
\be\la{evf}\ba
\begin{cases}
\Delta G=\div \left(\rho \dot{u} \right) &\mathrm{in}\, \,  \OM, \\
\frac {\p G}{\p n}=\left(\rho \dot{u}-\mu \na^{\bot} \o \right) \cdot n &\mathrm{on}\, \,  \p \OM.
\end{cases}
\ea\ee

Note that the Green's function $N(x,y)$ for the Neumann problem (see \cite{STT}) on the unit disc $\mathbb{D}$ is formulated as
\be\ba\nonumber
N(x,y)=-\frac{1}{2\pi}\bigg(\log|x-y|+\log\left||x|y-\frac{x}{|x|}\right|\bigg).
\ea\ee

Moreover, according to the Riemann mapping theorem (see \cite{SES}), there exists a conformal mapping
$\varphi=(\varphi_1, \varphi_2):\overline{\Omega}\rightarrow\overline{\mathbb{D}}$.
By using the Green's function on the unit disk and the conformal mapping, we can obtain the point-wise representation of $G$ on $\OM$.
We define the pull back Green's function $\widetilde{N}(x,y)$ of $\Omega$ as follows:
\be\ba\nonumber
\widetilde{N}(x,\, y)=N\big(\varphi(x),\varphi(y)\big) \ \ \mathrm{for}\ x,y\in\Omega.
\ea\ee

Employing the pull back Green's function $\widetilde{N}$, we may derive the point-wise representation of $G$, which is presented in the following lemma.
The proof can be found in \cite[Lemma 3.7]{FLL}.

\begin{lemma}\label{l7}
Assume that $G\in C\big([0,T];C^1(\bar\Omega)\cap C^2(\Omega)\big)$ satisfies the 
equation $(\ref{evf})$. Then for any $x\in\Omega$, it holds that
\be\la{nsb07}\ba
G(x,t)
=&-\int_\Omega \nabla_y\widetilde{N}(x,\, y) \cdot \n \du(y)\, dy
+\int_{\partial\Omega} \frac{\partial \widetilde{N}}{\partial n}(x,\, y) G(y)dS_y \\
&+\mu \int_{\partial\Omega} \widetilde{N}(x,y) n^\bot \cdot \nabla
\left( A u\cdot n^\bot \right) dS_y \\
=&-\frac{D}{Dt}\psi(x,t)+K_1(x,t)+K_2(x,t),
\ea\ee
where
\be\ba\nonumber
\psi(x,t)&\triangleq\int_\Omega\nabla_y\widetilde{N}\big(x,y\big)\cdot \rho u(y) dy,\\
K_1(x,t)& \triangleq \int_{\partial\Omega} \frac{\partial \widetilde{N}}{\partial n}(x,y) G(y) dS_y
+\mu \int_{\partial\Omega} \widetilde{N}(x,y) n^\bot \cdot \nabla
\left( A u\cdot n^\bot \right) dS_y, \\
K_2(x,t) & \triangleq\int_\Omega\left[ \p_{x_i} \p_{y_j}\widetilde{N}(x,y)\cdot u_i(x)+\p_{y_i} \p_{y_j}\widetilde{N}(x,y)\cdot u_i(y) \right]\rho u_j(y)dy.
\ea\ee
\end{lemma}

\begin{lemma}\label{l8}
Let $(\n,u)$ be a classical solution to \eqref{ns}--\eqref{i3} satisfying \eqref{p11}.
Then there exists a positive constant $C$ depending only on
$\ga$, $\mu$, $E_0$, $\| \n_0 \|_{L^\infty}$, $A$, $\OM$, and $\|\na u_0\|_{L^2}$ such that
\be\la{nsb08}\ba
\| \psi \|_{L^\infty} \le C \nu^{\frac{3}{4}},
\ea\ee
\be\la{nsb008}\ba
\| K_1 \|_{L^\infty} \le C \left( \frac{1}{\nu} + B_1 + B_2 \right),
\ea\ee
\be\la{nsb0008}\ba
\| K_2 \|_{L^\infty} 
\le C \nu^{\frac{1}{2}} \| \na u\|_{L^2}^{\frac{1}{2}}
\left( B_1^{\frac{1}{2}} B_2^{\frac{1}{2}} + B_1
+ \frac{1}{\nu} \right) + C B_1^2.
\ea\ee
\end{lemma}
\begin{proof}	
First, it follows from Lemma \ref{l3} that for $\de = \nu^{-\frac{1}{2}} \de_0$,
\be\ba\nonumber
\left|\int_\Omega\nabla_{y}\tilde{N}\cdot \rho u \, dy\right|
& \le C \int_\OM |x-y|^{-1} \n |u| dy \\
&\leq C\left(\int_\Omega|x-y|^{-\frac{2+\de}{1+\de}}dy\right)^{\frac{1+\de}{2+\de}}\left(\int_\Omega \rho^{2+\de}|u|^{2+\de}dy\right)^{\frac{1}{2+\de}}\\
&\leq C\de^{-\frac{1+\de}{2+\de}} \nu^{\frac{1}{2+\de}} \\
&\leq C\nu^{\frac{3}{4}},
\ea\ee
which yields (\ref{nsb08}).

For the first term of $K_1$, according to \cite[Lemma 3.6]{FLL},
it holds that for any $x\in \OM$, $y\in \p \OM$,
\be\ba\nonumber
\frac{\partial \widetilde{N}}{\partial n}(x,\, y)=-\frac{1}{2 \pi} |\na \varphi_1 (y)|,
\ea\ee
which together with (\ref{nsb45}) and Sobolev embedding implies that
\be\la{nsb84}\ba
\int_{\partial\Omega}\left|\frac{\partial \widetilde{N}}{\partial n}(x,\, y) G(y) \right|dS_y
\leq C\big\|G\big\|_{H^1}\leq C \left( B_1 + B_2 \right).
\ea\ee
On the other hand, we deduce from (\ref{i3}), (\ref{nsb04}) and H\"older's inequality that
\be\ba\nonumber
\left| \int _{\partial \Omega }{\widetilde{N}} n^\bot \cdot \nabla (Au\cdot n^\bot )\mathrm{d}S_y\right|
& =\left| \int _{\Omega }\mathrm {div}(\nabla ^\bot (Au\cdot n^\bot ){\widetilde{N}}) dy \right| \\
& \leq C\int _{\Omega }\vert \nabla {\widetilde{N}} \vert \cdot (\vert \nabla u\vert +\vert u\vert )\,\mathrm{d}x \\
& \leq C\Vert \nabla {\widetilde{N}}\Vert _{L^\frac{4}{3}} \Vert \nabla u\Vert _{L^4} \\
& \le C \left( B_1^{\frac{1}{2}} B_2^{\frac{1}{2}}
+ \frac{1}{\nu} \| P-\ol{P} \|_{L^4} + B_1 \right) \\
& \le C \left( B_1 + B_2 + \frac{1}{\nu} \right),
\ea\ee
which combined with (\ref{nsb84}) gives (\ref{nsb008}).

It remains to estimate $K_2$.
By \cite[Proposition 3.2]{FLL}, we have
\be\la{nsb87}\ba
\|K_2(\cdot, t)\|_{L^\infty(\Omega)}&\leq C\left(\sup_{x\in\overline{\Omega}}\int_\Omega\frac{\rho|u|^2(y)}{|x-y|}\, dy
+\sup_{x\in\overline{\Omega}}\int_\Omega\frac{|u(x)-u(y)|}{|x-y|^2} \cdot \rho|u|(y)\, dy\right).
\ea\ee

For the first term on the right side of (\ref{nsb87}), by applying 
H\"older's inequality and Sobolev embedding, we derive
\be\la{nsb88}\ba
\int_\Omega\frac{\rho|u|^2(y)}{|x-y|}dy
\leq C \left(\int_\Omega |x-y|^{-3/2}dy\right)^{2/3} \left(\int_\Omega |u|^{6} dy\right)^{1/3}
\leq C \|\nabla u\|_{L^2}^2.
\ea\ee

Moreover, by the Sobolev embedding theorem (Theorem 4 of \cite[Chapter 5]{EL}), we obtain that for any $x, y\in \overline{\Omega}$,
\bnn
|u(x)-u(y)| \leq C \|\nabla u\|_{L^4}|x-y|^{\frac{1}{2}},
\enn
which yields that
\be\la{nsb89}\ba
\int_\Omega\frac{|u(x)-u(y)|}{|x-y|^2}\rho|u|(y) dy
& \leq C \|\nabla u\|_{L^4} \int_\Omega |x-y|^{-\frac{3}{2}} \rho|u|(y)dy.
\ea\ee

For $\tau>0 $ and $\varepsilon_0\in (0, \frac{1}{4})$ to be determined later, we conclude from (\ref{gn11}) and H\"older's inequality that
\be\la{nsb810}\ba
& \int_{|x-y|<\tau}|x-y|^{-\frac{3}{2}} \rho|u|(y) dy \\
& \leq C \left(\int_{|x-y|<\tau} |x-y|^{-\frac{3}{2}(1+\varepsilon_0)} dy \right)^{\frac{1}{1+\varepsilon_0}}
\|u\|_{L^{\frac{1+\varepsilon_0}{\varepsilon_0}}} \\
&\leq C \varepsilon_0^{-\frac{1}{2}}
\tau^{\frac{1}{2}-\frac{2\varepsilon_0}{1+\varepsilon_0}} \|\na u\|_{L^2}.
\ea\ee
In addition, by using (\ref{nsb03}) and H\"older's inequality, we conclude that
\be\la{nsb811}\ba
&\int_{|x-y|>\tau}|x-y|^{-\frac{3}{2}} \rho|u|(y)dy \\
&\leq \left(\int_{|x-y|>\tau} |x-y|^{-\frac{3}{2} (\frac{2+\de}{1+\de})} dy\right)^{\frac{1+\de}{2+\de}}
\left(\int_\Omega\rho^{2+\de}|u|^{2+\de}dx\right)^{\frac{1}{2+\de}} \\
&\leq C \tau^{-\frac{1}{2}+\frac{\de}{2+\de}} \nu^{\frac{1}{2+\de}}.
\ea\ee
If $\| \na u \|_{L^2} = 0$, then $u=0$ by Poincar\'e's inequality, and the desired estimate is immediate.
We may therefore assume that $\| \na u \|_{L^2} > 0$ and choose $\tau>0$ such that
\be\la{nsb812}\ba
\tau^{-\frac{1}{2}+\frac{\de}{2+\de}} =
\| \na u\|_{L^2}^{ \frac{1}{2} },
\ea\ee
and set
\be\ba\nonumber
\varepsilon_0=\frac{\de}{4+\de} \in \left(0,\frac{1}{4}\right).
\ea\ee
Then
\be\la{nsb814}\ba 
\tau^{\frac{1}{2}-\frac{2\varepsilon_0}{1+\varepsilon_0}} \|\na u\|_{L^2} = \| \na u\|_{L^2}^{\frac{1}{2}}.
\ea\ee
Combining (\ref{nsb810}), (\ref{nsb811}), (\ref{nsb812}) and (\ref{nsb814}), we infer that
\be\la{nsb815}\ba
\int_\Omega |x-y|^{-\frac{3}{2}} \rho|u|(y) dy
& \leq C \varepsilon_0^{-\frac{1}{2}}  \| \na u\|_{L^2}^{\frac{1}{2}} + C \| \na u\|_{L^2}^{ \frac{1}{2} } \nu^{\frac{1}{2+\de}} \\
&\leq C \nu^{\frac{1}{2}} \| \na u\|_{L^2}^{ \frac{1}{2} },
\ea\ee
where in the last inequality we have used 
$\varepsilon_0^{-\frac{1}{2}}\le C(p)\de^{-\frac{1}{2}}\le C(p)\nu^{\frac{1}{4}}$
due to (\ref{nsb003}).

Therefore, it follows from (\ref{nsb89}), (\ref{nsb815}) and (\ref{nsb04}) that
\be\la{nsb816}\ba
& \int_\Omega\frac{|u(x )-u(y)|}{|x-y|^2}\rho|u|(y)dy \\
& \le C \nu^{\frac{1}{2}} \|\nabla u\|_{L^4} \| \na u\|_{L^2}^{\frac{1}{2}} \\
&\leq C \nu^{\frac{1}{2}} \| \na u\|_{L^2}^{\frac{1}{2}} \left( B_1^{\frac{1}{2}} B_2^{\frac{1}{2}}
+ B_1 + \frac{1}{\nu} \| P-\ol{P} \|_{L^4} \right),
\ea\ee
which together with (\ref{nsb88}) implies (\ref{nsb0008}) and completes the proof of Lemma \ref{l8}.
\end{proof}

With Lemmas \ref{l1}--\ref{l8} at hand, we now proceed to the proof of Proposition \ref{pl1}.
\begin{proof}[Proof of Proposition \ref{pl1}]
First, from $(\ref{ns})_1$ and (\ref{gw}), we obtain
\be\ba\la{nsb91}
\pa_t \log\n + u\cdot\na \log\n + \frac{1}{\nu}( P-\ol{P} + G) = 0.
\ea\ee
Inserting (\ref{nsb07}) into (\ref{nsb91}) yields
\be\ba\nonumber
\pa_t F + u\cdot\na F + \frac{1}{\nu}(P-\ol{P}) =
-\frac{K_1}{\nu} -\frac{K_2}{\nu},
\ea\ee
with
\be\ba\la{F}
F \triangleq \log\n - \frac{1}{\nu}\psi.
\ea\ee
Note that $\n^\ga \geq \ga\log\n + 1$, which implies that
\be\ba\la{nsb93}
\pa_t F + u\cdot\na F + \frac{\ga}{\nu} F
\leq \frac{\ga}{\nu^2}|\psi| + \frac{1}{\nu}\left \lvert K_1 \right\rvert 
+\frac{1}{\nu} \left\lvert K_2 \right\rvert
+ \frac{1}{\nu |\om|} \ol{P}.
\ea\ee

In order to handle the material derivative $\frac{D}{Dt}F =\pa_t F + u\cdot \na F$,
we introduce the characteristic curve $y(s;x,t)$ which is defined by
\be\ba\nonumber
\begin{cases}
	\frac{d}{ds}y(s) = u(y,s),\\
	y(t;x,t) = x,
\end{cases}
\ea\ee
which together with (\ref{nsb93}) yields for all $(x,t)\in\OM \times (0,T]$,
\be\la{nsb94}\ba
\frac{d}{ds}F(s) + \frac{\ga}{\nu} F(s)
\leq \frac{\ga}{\nu^2} \|\psi\|_{L^\infty} + \frac{1}{\nu}\| K_1 \|_{L^\infty}+\frac{1}{\nu}\| K_2 \|_{L^{\infty}}
+ \frac{\ga-1}{\nu |\om|} E_0,
\ea\ee
where we denote slightly abusively $F(y(s;x,t),s)$ by $F(s)$.

Using (\ref{nsb008}), (\ref{nsb0008}), and Young's inequality, we derive
\be\la{nsb95}\ba
\frac{1}{\nu} \| K_1 \|_{L^{\infty}}
\le \frac{C}{\nu} \left( \frac{1}{\nu} + B_1 + B_2 \right)
\le C \nu^{-\frac{3}{2}} + C \nu^{-\frac{1}{2}} B^2_1 + \frac{C}{\nu} B_2,
\ea\ee
and
\be\la{nsb950}\ba
\frac{1}{\nu} \| K_2 \|_{L^{\infty}}
& \le C \nu^{-\frac{1}{2}} B^{\frac{1}{2}}_1 \left( B^{\frac{1}{2}}_1 B^{\frac{1}{2}}_2
+ B_1 + \frac{1}{\nu} \right) + C \nu^{-\frac{1}{2}} B^2_1 \\
& \le C \nu^{-\frac{1}{4}} B^2_1 + C \nu^{-\frac{3}{4}} B_2 + C \nu^{-\frac{5}{4}}.
\ea\ee
Combining (\ref{nsb08}), (\ref{nsb94}), (\ref{nsb95}) and (\ref{nsb950}) yields that
\be\ba\nonumber
\frac{d}{ds}F(s) + \frac{\ga}{\nu} F(s)
\le C \nu^{-\frac{5}{4}} + C \nu^{-\frac{1}{4}} B^2_1 + C \nu^{-\frac{3}{4}} B_2 
+ \frac{\ga-1}{\nu |\om|} E_0.
\ea\ee
An application of the maximum principle shows that
\be\ba\la{nsb96}
F(t) \le & e^{-\frac{\ga}{\nu}t}F(0) +C \nu^{-\frac{5}{4}}\int_0^t e^{-\frac{\ga}{\nu}(t-s)} ds
+ C \nu^{-\frac{1}{4}} \int_0^t e^{-\frac{\ga}{\nu}(t-s)} B^2_1 ds \\
& + C \nu^{-\frac{3}{4}} \int_0^t e^{-\frac{\ga}{\nu}(t-s)} B_2 ds
+\frac{\ga-1}{\ga |\om|}(1-e^{-\frac{\ga}{\nu}t}) E_0 \\
\le & e^{-\frac{\ga}{\nu}t}F(0) + C \nu^{-\frac{1}{6}} + \frac{\ga-1}{\ga |\om|} E_0,
\ea\ee
where in the second inequality we have used the following estimate:
\be\ba\nonumber
\int_0^t e^{-\frac{\ga}{\nu}(t-s)} B_2 ds
& = \int_0^t e^{-\frac{\ga}{\nu}(t-s)}
\frac{ B_2 }{(2+B^2_1)^{\frac{1}{2}}} (2+B^2_1)^{\frac{1}{2}} ds \\
& \le \left(\int_0^t e^{-\frac{2 \ga}{\nu}(t-s)} (2 + B^2_1) ds \right)^{\frac{1}{2}}
\left(\int_0^t \frac{B_2^2}{2+B^2_1} ds \right)^{\frac{1}{2}} \\
& \le C \nu^{\frac{1}{2}} \log^{\frac{1}{2}}(2+\nu) \\
& \le C \nu^{\frac{7}{12}},
\ea\ee
due to (\ref{nsb06}) and H\"older's inequality.

Consequently, we deduce from (\ref{F}), (\ref{nsb96}) and (\ref{nsb08}) that
\be\ba\la{nsb920}
\log \n &= F+\frac{1}{\nu} \psi \\
& \le e^{-\frac{\ga}{\nu}t} \log \|\n_0\|_{L^\infty} + \frac{1}{\nu} \left( \| \psi_0 \|_{L^\infty}+\| \psi \|_{L^\infty} \right)
+ C \nu^{-\frac{1}{6}} + \frac{\ga-1}{\ga |\om|} E_0 \\
& \le e^{-\frac{\ga}{\nu}t} \log \|\n_0\|_{L^\infty} + C \nu^{-\frac{1}{6}}
+\frac{\ga-1}{\ga |\om|} E_0.
\ea\ee
Since
\be\ba\nonumber
\|\n_0\|_{L^\infty} \ge \int \n_0 dx =1,
\ea\ee
we conclude that
\be\ba\nonumber
\log \n \le \log \|\n_0\|_{L^\infty} + C_0 \nu^{-\frac{1}{6}} + \frac{\ga-1}{\ga |\om|} E_0,
\ea\ee
where the constant $C_0$ depends only on
$\ga,\ \mu,\ A,\ \OM,\ \|\n_0\|_{L^\infty},\ E_0$ and $\|\na u_0\|_{L^2}$, 
but is independent of $T$ and $\nu$.
Finally, we choose
\be\ba\la{nsb923}
\nu_1 = \left( \frac{ C_0 }{\log \frac{3}{2}} \right)^6, 
\ea\ee
then when $\nu\geq\nu_1$, we obtain (\ref{p11}).
\end{proof}

\section{A Priori Estimates (II): Higher Order Estimates}
In this section, we assume $\nu \ge \nu_1$, with $\nu_1$ defined in Proposition \ref{pl1}.
Let $(\n,u)$ be a smooth solution of (\ref{ns})--(\ref{i3}) on $\OM \times (0,T]$ satisfying (\ref{p11}).
We establish higher-order estimates and exponential decay for $(\n,u)$, following arguments similar to those in \cite{CL,HLX2,LX3,LLL}.
\begin{lemma}\la{g1}
There exists a positive constant $C$ depending only on
$\mu$, $\gamma$, $\OM$, $A$, $E_0$, $\| \n_0 \|_{L^\infty}$, and $\| u_0 \|_{H^1}$ such that
\be\ba\la{gb01}
\sup_{0\le t\le T}
 \si \int\n|\dot u|^2dx
+\int_0^{T} \si \| \na\dot u\|^2_{L^2} dt \le C\nu^C,
\ea\ee
and, for any $\tau \in (0,T)$,
\be\ba\la{gb01a}
\sup_{\tau \le t \le T} \int \n |\dot u|^2 dx
+ \int_\tau^{T} \| \na\dot u\|^2_{L^2} dt \le C (1 + \tau^{-C}) e^{C\tau^{-C}}.
\ea\ee
Moreover, if the initial data $(\n_0,u_0)$ satisfy 
\be\la{gb001}\ba
\div u_0 = 0,
\ea\ee 
then
\be\la{gb002}\ba
\sup_{0\le t\le T}
\si \int\n|\dot u|^2dx
+\int_0^{T} \si \| \na\dot u\|^2_{L^2} dt \le C.
\ea\ee
\end{lemma}
\begin{proof}
The idea of this proof comes from \cite{CL,H1,FLL,LX3}.
Operating $ \dot u^j[\frac{\pa}{\pa t}+\div(u\cdot)]$ on $ (\ref{ns})_2^j,$ summing with respect to $j,$
and integrating the resulting equation over $\OM$, we obtain after integration by parts that
\be\la{gb11}\ba
\frac{d}{dt}\left(\frac{1}{2}\int\rho|\dot{u}|^2dx \right)
&=\int \bigg( {\dot{u}}\cdot \nabla G_t +{\dot{u}}^j\mathrm {div}(u \partial _jG) \bigg) dx\\
&\quad +\mu \int \bigg( {\dot{u}}\cdot \nabla ^\bot \omega _t 
+{\dot{u}}^j\partial _k(u^k (\nabla ^\bot \omega )_j ) \bigg) dx \\
&=I_1+I_2.
\ea\ee
For $I_1$, integration by parts with H\"older's and Young inequalities yields
\be\la{gb12}\ba
I_1 & = \int_{\partial \Omega} G_t ( {\dot{u}}\cdot n) ds
- \int \div \dot{u} \left( \dot{G}-u \cdot \na G \right) dx
+ \int u \cdot \na \dot{u}^j \p_j G dx \\
& \le \int_{\partial \Omega} G_t ( {\dot{u}}\cdot n) ds
- \int \div \dot{u} \dot{G} dx + C \| \na \dot{u} \|_{L^2} \| u\|_{L^6} \| \na G \|_{L^3} \\
& \le  \int_{\partial \Omega} G_t ( {\dot{u}}\cdot n) ds
- \int \div \dot{u} \dot{G} dx
+ \varepsilon \Vert \nabla {\dot{u}}\Vert _{L^2}^2 \\
& \quad +C(\varepsilon) \Vert \sqrt{\rho }{\dot{u}}\Vert _{L^2}^2 \left(\| \na u\|^4_{L^2}+\| \na u\|^2_{L^2} \right) 
+ C(\varepsilon) \left( \| \na u\|^2_{L^2}+\| \na u\|^6_{L^2} \right),
\ea\ee
where in the last inequality we have used the following estimate:
\be\la{gb14}\ba
\| \na G \|_{L^3}+\| \na \o \|_{L^3} & \le \| \na G \|^{\frac{1}{2}}_{L^2} \| \na G \|^{\frac{1}{2}}_{L^6}+\| \na \o \|^{\frac{1}{2}}_{L^2} \| \na \o \|^{\frac{1}{2}}_{L^6} \\
& \le C\left(\| \sqrt{\n} \dot{u}\|_{L^2} + \| \na u \|_{L^2} \right)^{\frac{1}{2}} 
\left(\| \n \dot{u}\|_{L^6} + \| \na u \|_{L^6} \right)^{\frac{1}{2}} \\
& \le C\left(\| \sqrt{\n} \dot{u}\|_{L^2} + \| \na u \|_{L^2} \right)^{\frac{1}{2}} 
\left(\| \na \dot{u}\|_{L^2}+\| \na u\|^2_{L^2} \right)^{\frac{1}{2}} \\
& \quad +C\left(\| \sqrt{\n} \dot{u}\|_{L^2} + \| \na u \|_{L^2} \right)^{\frac{1}{2}} 
\left( 1+\| \sqrt{\n} \dot{u}\|_{L^2} + \| \na u \|_{L^2} \right)^{\frac{1}{2}} \\
& \le C\left(\| \sqrt{\n} \dot{u}\|_{L^2} + \| \na u \|_{L^2} \right)^{\frac{1}{2}} 
\left(\| \na \dot{u}\|_{L^2}+\| \na u\|^2_{L^2} + 1 \right)^{\frac{1}{2}},
\ea\ee
due to (\ref{nsb44}), (\ref{nsb45}), (\ref{dc1}), (\ref{emdu1}) and H\"older's inequality.

Next, for the boundary term of (\ref{gb12}), with the help of (\ref{i3}) and (\ref{bjds}), we derive
\be\la{gb15}\ba
&\int _{\partial \Omega }G_t({\dot{u}}\cdot n) ds\\ 
&= \int _{\partial \Omega }G_t(u\cdot \nabla u\cdot n) ds \\
&= -\frac{d}{dt} \int _{\partial \Omega } G (u\cdot \nabla n\cdot u)ds
+\int _{\partial \Omega } G(u\cdot \nabla n\cdot u)_t ds\\ 
&= -\frac{d}{dt} \int _{\partial \Omega } G(u\cdot \nabla n\cdot u)ds
+\int _{\partial \Omega } G(u_t\cdot \nabla n\cdot u) ds
+\int _{\partial \Omega } G(u\cdot \nabla n\cdot u_t)ds \\ 
&= -\frac{d}{dt} \int _{\partial \Omega } G(u\cdot \nabla n\cdot u)ds
+ \int _{\partial \Omega } G({\dot{u}}\cdot \nabla n\cdot u) + G(u\cdot \nabla n\cdot {\dot{u}}) ds  \\
&\quad - \int _{\partial \Omega } G \left( (u\cdot \nabla u)\cdot \nabla n\cdot u\right) ds-\int _{\partial \Omega } G \left( u\cdot \nabla n\cdot (u\cdot \nabla u) \right) ds \\
&=-\frac{d}{dt} \int_{\partial \Omega } G (u\cdot \nabla n\cdot u)ds+J_1+J_2+J_3.
\ea\ee
For $J_1$, it follows from (\ref{emdu1}), (\ref{nsb45}) and Young's inequality that
\be\la{gb16}\ba
J_1&=\int _{\partial \Omega } G({\dot{u}}\cdot \nabla n\cdot u) + G(u\cdot \nabla n\cdot {\dot{u}}) ds \\
&\leq C\Vert u\Vert _{H^1}\Vert {\dot{u}}\Vert _{H^1}\Vert G \Vert _{H^1}\\
&\leq C \Vert \nabla u\Vert _{L^2} \left( \Vert \nabla {\dot{u}}\Vert _{L^2}+\Vert \nabla u\Vert _{L^2}^2 \right) 
\left( \Vert \sqrt{\rho }{\dot{u}}\Vert _{L^2}+\Vert \nabla u\Vert _{L^2} \right) \\
&\leq \varepsilon \Vert \nabla {\dot{u}}\Vert _{L^2}^2 +C(\varepsilon ) \Vert \sqrt{\rho }{\dot{u}}\Vert _{L^2}^2 \| \na u\|^2_{L^2} + C(\varepsilon )\| \na u\|^4_{L^2}.
\ea\ee
To deal with $J_2$, by virtue of (\ref{nsb45}), (\ref{bjds}) and Young's inequality, we arrive at
\be\la{gb17}\ba
\vert J_2\vert= & \left| \int _{\partial \Omega } G((u\cdot \nabla u)\cdot \nabla n\cdot u) ds \right| \\
= & \left| \int _{\partial \Omega } G (u\cdot n^\bot )n^\bot \cdot \nabla u^i\partial_i n_j u^j ds \right|  \\
= & \left| \int \nabla ^\bot \cdot \left( \nabla u^i \partial_i n_j u^j  G (u\cdot n^\bot )\right) dx \right| \\
= & \left| \int \nabla u_i\cdot \nabla ^\bot \left( \partial_i n_j u^j G (u\cdot n^\bot )\right) dx \right| \\
\le & C\int \vert \nabla u\vert \left( \vert G \vert \vert u\vert ^2+\vert G\vert \vert u\vert \vert \nabla u\vert +\vert u\vert ^2\vert \nabla G \vert \right) dx \\
\le & C\Vert \nabla u \Vert _{L^4}\left( \Vert G \Vert _{L^4}\Vert u\Vert _{L^4}^2+\Vert G \Vert _{L^4}\Vert u\Vert _{L^4}\Vert \nabla u\Vert _{L^4}+\Vert \nabla G \Vert _{L^2}\Vert u\Vert _{L^8}^2\right)  \\
\le & C \left(\Vert \nabla u \Vert _{L^4} \Vert \nabla u \Vert^2_{L^2}+\Vert \nabla u \Vert^2_{L^4} \Vert \nabla u \Vert_{L^2} \right)  \Vert G \Vert _{H^1} \\
\le & C\Vert \nabla u \Vert^2_{L^4} \Vert \nabla u \Vert_{L^2}
\left(\Vert \sqrt{\rho }{\dot{u}}\Vert_{L^2} + \| \na u\|_{L^2} \right) \\
\le & C \Vert \sqrt{\rho }{\dot{u}}\Vert _{L^2}^2 \| \na u\|^2_{L^2}+C \Vert \nabla u \Vert^4_{L^4}.
\ea\ee
Similarly, we also have
\be\la{gb18}\ba
\vert J_3\vert \le C \Vert \sqrt{\rho }{\dot{u}}\Vert _{L^2}^2 \| \na u\|^2_{L^2}+C \Vert \nabla u \Vert^4_{L^4},
\ea\ee
which together with (\ref{gb15})--(\ref{gb17}) implies that
\be\la{gb19}\ba
\int _{\partial \Omega }G_t({\dot{u}}\cdot n) ds 
\leq &-\frac{d}{dt}\int _{\partial \Omega } G (u\cdot \nabla n\cdot u)ds
+\varepsilon \Vert \nabla {\dot{u}}\Vert _{L^2}^2\\
&+C(\varepsilon ) \Vert \sqrt{\rho }{\dot{u}}\Vert_{L^2}^2 \| \na u\|^2_{L^2} + C(\varepsilon )\| \na u\|^4_{L^4}.
\ea\ee
For the first term in the last line of (\ref{gb12}), we deduce from (\ref{gw})
and $(\ref{ns})_1$ that
\be\ba\nonumber
\div \dot{u} & = (\div u)_t + \p_i u^j \p_j u^i + u \cdot \na \div u \\
& = \frac{1}{\nu} \dot{G} + \frac{1}{\nu} (P_t+u \cdot \na P) - \frac{1}{\nu} \ol{P_t}
+ \p_i u^j \p_j u^i \\
& = \frac{1}{\nu} \dot{G} - \frac{\ga}{\nu} P \div u + \frac{\ga-1}{\nu} \int P \div u dx + \p_i u^j \p_j u^i,
\ea\ee
which together with Young's inequality yields
\be\la{gb109}\ba
- \int \div \dot{u} \dot{G} dx
& = -\frac{1}{\nu} \| \dot{G} \|^2_{L^2} 
- \int \dot{G} \p_i u^j \p_j u^i dx \\
& \quad + \frac{\ga}{\nu} \int \dot{G} P \div u dx
-\frac{\ga-1}{\nu} \int P \div u dx \int \dot{G} dx \\
& \le -\frac{1}{2\nu} \| \dot{G} \|^2_{L^2}+\frac{C}{\nu} \| \div u \|^2_{L^2}
- \int \dot{G} \p_i u^j \p_j u^i dx.
\ea\ee
In addition, by integrating by parts and using H\"older's inequality, we obtain
\be\la{gb1009}\ba
-\int \dot{G} \p_i u^j \p_j u^i dx 
& = -\int ( G_t + u \cdot \na G ) \p_i u \cdot \na u^i dx \\
& = -\frac{d}{dt} \left( \int G \p_i u \cdot \na u^i dx \right) 
+ 2 \int G \p_i u \cdot \na u^i_t dx \\
& \quad + \int G \div u \p_i u \cdot \na u^i dx 
+ 2 \int G u \cdot \na \p_i u \cdot \na u^i dx \\
& = -\frac{d}{dt} \left( \int G \p_i u \cdot \na u^i dx \right)
+ 2 \int G \p_i u \cdot \na \dot{u}^i dx \\
& \quad - 2 \int G \p_i u \cdot \na u \cdot \na u^i dx
+ \int G \div u \p_i u \cdot \na u^i dx \\
& \le -\frac{d}{dt} \left( \int G \p_i u \cdot \na u^i dx \right)
+ C \| \na \dot{u} \|_{L^2} \| \na u \|_{L^4} \| G \|_{L^4}
+ C \| \na u \|^3_{L^4} \| G \|_{L^4} \\
& \le -\frac{d}{dt} \left( \int G \p_i u \cdot \na u^i dx \right)
+ \ep \| \na \dot{u} \|^2_{L^2}+C(\ep) \| \na u \|^4_{L^4}+C(\ep) \| G \|^4_{L^4}.
\ea\ee
Combining (\ref{gb12}), (\ref{gb109}) and (\ref{gb1009}) implies
\be\la{gb110}\ba
I_1 & \le -\frac{d}{dt} \left( \int G \p_i u \cdot \na u^i dx + \int _{\partial \Omega } G (u\cdot \nabla n\cdot u)ds \right)
-\frac{1}{2\nu} \| \dot{G} \|^2_{L^2} +3\varepsilon \Vert \nabla {\dot{u}}\Vert _{L^2}^2 \\
& \quad + C(\ep) \| \na u \|^4_{L^4}
+ C(\ep) \| G \|^4_{L^4} + C(\varepsilon) \left( \| \na u\|^2_{L^2}+\| \na u\|^6_{L^2} \right) \\
& \quad + C(\varepsilon) \Vert \sqrt{\rho }{\dot{u}}\Vert _{L^2}^2 
\left(\| \na u\|^4_{L^2} + \| \na u\|^2_{L^2} \right).
\ea\ee

For $I_2$, integration by parts and applying (\ref{gb14}) gives
\be\la{gb111}\ba
I_2&=\mu \int \left({\dot{u}}\cdot \nabla ^\bot \omega _t
+{\dot{u}}^j \p_k \left( u^k (\na^\bot \o)^j \right) \right) dx \\
&=\mu \int _{\partial \Omega }({\dot{u}}\cdot n^\bot )\omega _t ds
-\mu \int \curl \dot{u} \o_t dx-\mu \int u \cdot \na \dot{u} \cdot \na^\bot \o dx \\
&=\mu \int_{\p \OM} -A ({\dot{u}}\cdot n^\bot)^2 + ({\dot{u}}\cdot n^\bot) A (u \cdot \na u \cdot n^\bot) ds
-\mu \int (\curl \dot{u})^2 dx \\
&\quad + \mu \int \curl \dot{u} \curl (u \cdot \na u) dx
-\mu \int u \cdot \na \dot{u} \cdot \na^\bot \o dx \\
&\le \mu \int_{\p \OM} -A ({\dot{u}}\cdot n^\bot)^2 + ({\dot{u}}\cdot n^\bot) A (u \cdot \na u \cdot n^\bot) ds
-\mu \int (\curl \dot{u})^2 dx \\
&\quad + C \int |\na \dot{u}| |\na u|^2 + |\na \dot{u}| |\na \o| |u| dx \\
& \le -\mu \| \curl \dot{u}\|^2_{L^2} + 2\ep \| \na \dot{u} \|^2_{L^2} + C(\ep) \| \na u \|^4_{L^4}
+ C \| \na \dot{u} \|_{L^2} \| \na \o \|_{L^3} \| u \|_{L^6} \\
& \le -\mu \| \curl \dot{u}\|^2_{L^2} + 3\ep \| \na \dot{u} \|^2_{L^2} + C(\ep) \| \na u \|^4_{L^4}
+ C(\varepsilon) \left( \| \na u\|^2_{L^2}+\| \na u\|^6_{L^2} \right) \\
& \quad + C(\varepsilon) \Vert \sqrt{\rho }{\dot{u}}\Vert _{L^2}^2 \left(\| \na u\|^4_{L^2}+\| \na u\|^2_{L^2} \right),
\ea\ee
where we have used the following fact:
\be\la{gb112}\ba
\mu \int_{\p \OM} ({\dot{u}}\cdot n^\bot) A (u \cdot \na u \cdot n^\bot) ds
&= \mu \int_{\p \OM} ({\dot{u}}\cdot n^\bot) A ( (u \cdot n^\bot) n^\bot \cdot \na u \cdot n^\bot) ds \\
& = \mu \int \na^\bot \cdot \left( \na u \cdot n^\bot ({\dot{u}}\cdot n^\bot) A (u \cdot n^\bot) \right) dx \\
& \le C \int |\na u| |\dot{u}| |u| + |\na \dot{u}| |\na u| |u| + |\na u|^2 |\dot{u}| dx \\
& \le C \| \dot{u} \|_{L^2} \| \na u \|^2_{L^4} + C \| \na \dot{u} \|_{L^2} \| \na u \|^2_{L^4} \\
& \le \ep \| \na \dot{u} \|^2_{L^2} + C(\ep) \| \na u \|^4_{L^4},
\ea\ee
due to (\ref{bjds}), (\ref{emdu1}) and Young's inequality.

Combining with (\ref{gb11}), (\ref{gb110}) and (\ref{gb111}), we arrive at
\be\la{gb116}\ba
& \frac{d}{dt} \left(\int\rho|\dot{u}|^2 dx + 2 \int G \p_i u \cdot \na u^i dx
+ 2 \int _{\partial \Omega } G (u\cdot \nabla n\cdot u)ds \right)
+ \frac{1}{\nu} \| \dot{G} \|^2_{L^2} + \mu \| \curl \dot{u}\|^2_{L^2} \\
& \le 12 \varepsilon \Vert \nabla {\dot{u}}\Vert _{L^2}^2
+ C(\ep) \| \na u \|^4_{L^4}
+ C(\ep) \| G \|^4_{L^4} + C(\varepsilon) \left( \| \na u\|^2_{L^2}+\| \na u\|^6_{L^2} \right) \\
& \quad + C(\varepsilon) \Vert \sqrt{\rho }{\dot{u}}\Vert _{L^2}^2 
\left(\| \na u\|^4_{L^2} + \| \na u\|^2_{L^2} \right) \\
& \le C \varepsilon \left( \| \div \dot{u}\|^2_{L^2}+\| \curl \dot{u}\|^2_{L^2} \right)
+ C(\ep) \| \na u \|^4_{L^4} + C(\ep) \| G \|^4_{L^4} \\
& \quad +C(\varepsilon) \left( \| \na u\|^2_{L^2}+\| \na u\|^6_{L^2} \right)
+ C(\varepsilon) \Vert \sqrt{\rho }{\dot{u}}\Vert _{L^2}^2 \left(\| \na u\|^4_{L^2}+\| \na u\|^2_{L^2} \right) \\
& \le C \varepsilon \left( \frac{1}{\nu^2} \| \dot{G} \|^2_{L^2} + \| \curl \dot{u}\|^2_{L^2} \right)
+ C(\ep) \| \na u \|^4_{L^4} + C(\ep) \| G \|^4_{L^4} \\
& \quad +C(\varepsilon) \left( \| \na u\|^2_{L^2}+\| \na u\|^6_{L^2} \right)
+ C(\varepsilon) \Vert \sqrt{\rho }{\dot{u}}\Vert _{L^2}^2 \left(\| \na u\|^4_{L^2}+\| \na u\|^2_{L^2} \right),
\ea\ee
where in the last inequality we have used (\ref{emdu2}) and the following fact:
\be\la{gb1106}\ba
\| \div \dot{u} \|^2_{L^2} \le \frac{2}{\nu^2} \| \dot{G} \|^2_{L^2} 
+ C \| \na u \|^2_{L^2} + C \| \na u \|^4_{L^4}
\ea\ee
Set
\be\la{gb11006}\ba
f(t) \triangleq \int\rho|\dot{u}|^2 dx + 2 \int G \p_i u \cdot \na u^i dx
+ 2 \int _{\partial \Omega } G (u\cdot \nabla n\cdot u)ds.
\ea\ee
By (\ref{nsb04}), we have
\be\la{gb116a}\ba
\| \na u \|_{L^4}^4
& \le C B_1^2 A_2^2 + C B_1^4 + \frac{C}{\nu^4} \| P - \ol{P} \|_{L^4}^4 \\
& \le C \| \sqrt{\n} \dot{u} \|_{L^2}^2 \left( 1 + \| \na u \|_{L^2}^2 + \nu \| \div u \|_{L^2}^2 \right) \\
& \quad + C \| \na u \|_{L^2}^4 + C \nu^2 \| \div u \|_{L^2}^4
+ \frac{C}{\nu} \| P - \ol{P} \|_{L^2}^2.
\ea\ee
Choosing $\varepsilon$ sufficiently small in (\ref{gb116}) and using (\ref{gb11006}) and (\ref{gb116a}), one obtains
\be\la{gb116b}\ba
& \frac{d}{dt} f(t)
+ \frac{1}{2\nu} \| \dot{G} \|^2_{L^2} + \frac{\mu}{2} \| \curl \dot{u}\|^2_{L^2} \\
& \le C \| G \|^4_{L^4}
+ C \left( \| \na u\|^2_{L^2} + \| \na u\|^6_{L^2} + \nu^2 \| \div u \|_{L^2}^4 \right) \\
& \quad + C \| \sqrt{\rho}{\dot{u}} \|_{L^2}^2 \left( 1 + \| \na u\|^4_{L^2} + \nu \| \div u\|^2_{L^2} \right)
+ \frac{C}{\nu} \| P - \ol{P} \|_{L^2}^2.
\ea\ee
Multiplying (\ref{gb116b}) by $\si$ leads to
\be\la{gb117}\ba
& \frac{d}{dt}\left( \si f(t) \right)
+\frac{1}{4} \si \left( \frac{1}{\nu} \| \dot{G} \|^2_{L^2} + \mu \| \curl \dot{u}\|^2_{L^2} \right) \\
& \le {\si}' f(t) + C \si \| G \|^4_{L^4}
+ C \si \left( \| \na u\|^2_{L^2} + \| \na u\|^6_{L^2} + \nu^2 \| \div u \|_{L^2}^4 \right) \\
& \quad + C \si \| \sqrt{\rho}{\dot{u}} \|_{L^2}^2 \left( 1 + \| \na u\|^4_{L^2} + \nu \| \div u\|^2_{L^2} \right)
+ \frac{C}{\nu} \| P - \ol{P} \|_{L^2}^2.
\ea\ee
Furthermore, we conclude from (\ref{nsb04}), (\ref{nsb45}) and Young's inequality that
\be\la{gb118}\ba
& \left| \int G \p_i u \cdot \na u^i dx + \int _{\partial \Omega } G (u\cdot \nabla n\cdot u) ds \right| \\
& \le C \| G \|_{L^4} \| \na u \|^2_{L^{\frac{8}{3}}} + C \| G\|_{H^1} \| \na u\|^2_{L^2} \\
& \le C \| G \|_{H^1} \left( B^{\frac{3}{2}}_1 B^{\frac{1}{2}}_2 + B^2_1 + \| P-\ol{P} \|^2_{L^{\frac{8}{3}}} \right) \\
& \le C \left(\| \sqrt{\n} \dot{u}\|_{L^2}+\| \na u\|_{L^2} \right)
\left( B^{\frac{3}{2}}_1 B^{\frac{1}{2}}_2 + B^2_1 + 1 \right) \\
& \le \frac{1}{2} \| \sqrt{\n} \dot{u}\|^2_{L^2}
+ C \left( B^6_1 + 1 \right).
\ea\ee
Integrating (\ref{gb117}) over $(0,T)$ and using (\ref{nsb1}), (\ref{nsb02}), (\ref{nsb06}), and (\ref{gb118}), we derive
\be\la{gb119}\ba
& \sup_{0 \le t \le T} \si \int \rho|\dot{u}|^2dx + \int_0^T \si 
\left( \frac{1}{\nu} \| \dot{G} \|^2_{L^2} + \mu \| \curl \dot{u}\|^2_{L^2} \right) dt 
\le C\nu^C.
\ea\ee
In addition, by virtue of (\ref{emdu2}), (\ref{nsb04}), and (\ref{gb119}), we arrive at
\be\la{gb120}\ba
\int_0^T \si \| \na \du \|^2_{L^2}
& \le C \int_0^T \si \left( \| \div \du \|^2_{L^2} + \| \curl \du \|^2_{L^2} 
+ \| \na u \|^4_{L^4} \right) dt \\
& \le C \int_0^T \si \left( \frac{1}{\nu^2} \| \dot{G} \|^2_{L^2} + \| \curl \du \|^2_{L^2} 
+ \| \na u \|^4_{L^4} + \| \na u \|^2_{L^2} \right) dt \\
& \le C\nu^C,
\ea\ee
which together with (\ref{gb119}) gives (\ref{gb01}).

Furthermore, by (\ref{gb11006}), (\ref{gb118}), and (\ref{nsb006}), we have for any $t \in (\frac{\tau}{2},T)$,
\be\la{ngb121}\ba
f(t) \ge \frac{1}{2} \| \sqrt{\n} \dot{u}(\cdot,t) \|_{L^2}^2 - C - C \tau^{-C}.
\ea\ee
For any $\tau \in (0,T)$, it follows from (\ref{gb116b}) , and Poincar\'e's inequality that for any $t \in (\frac{\tau}{2},T)$,
\be\la{ngb122}\ba
& \frac{d}{dt} f(t)
+ \frac{1}{2\nu} \| \dot{G} \|^2_{L^2} + \frac{\mu}{2} \| \curl \dot{u}\|^2_{L^2} \\
& \le C \| \sqrt{\rho}{\dot{u}} \|_{L^2}^4
+ C \left( \| \na u\|^2_{L^2} + \| \na u\|^6_{L^2} + \nu^2 \| \div u \|_{L^2}^4 \right) \\
& \quad + C \| \sqrt{\rho}{\dot{u}} \|_{L^2}^2 \left( 1 + \| \na u\|^4_{L^2} + \nu \| \div u\|^2_{L^2} \right)
+ \frac{C}{\nu} \| P - \ol{P} \|_{L^2}^2 \\
& \le C f(t) \| \sqrt{\rho}{\dot{u}} \|_{L^2}^2
+ C (1 + \tau^{-C}) \left( \| \na u\|^2_{L^2} + \nu \| \div u \|_{L^2}^2 \right) \\
& \quad + C (1 + \tau^{-C}) \| \sqrt{\rho}{\dot{u}} \|_{L^2}^2
+ \frac{C}{\nu} \| P - \ol{P} \|_{L^2}^2.
\ea\ee
By (\ref{nsb006}), there exists $t_* \in (\frac{\tau}{2},\tau)$ such that
\be\la{ngb123}\ba
\| \sqrt{\n} \dot{u}(\cdot,t_*) \|_{L^2}^2 \le \frac{2}{\tau} \int_{\frac{\tau}{2}}^{\tau} \| \sqrt{\n} \dot{u}(\cdot,t) \|_{L^2}^2 dt
\le \frac{2}{\tau} \int_{\frac{\tau}{2}}^{T} \| \sqrt{\n} \dot{u}(\cdot,t) \|_{L^2}^2 dt
\le C + C \tau^{-C},
\ea\ee
which together with (\ref{gb11006}), (\ref{gb118}), (\ref{nsb006}), and (\ref{ngb123}) yields
\be\la{ngb124}\ba
|f(t_*)| \le C + C \tau^{-C}.
\ea\ee
Applying Gr\"onwall's inequality to (\ref{ngb122}) over $(t_*,T)$ and using (\ref{nsb006}) and (\ref{nsb02}), we derive
\be\la{ngb125}\ba
& \sup_{t_* \le t \le T} f(t) + \int_{t_*}^T \left( \frac{1}{2\nu} \| \dot{G} \|^2_{L^2} + \frac{\mu}{2} \| \curl \dot{u}\|^2_{L^2} \right) dt \\
& \le C (1 + \tau^{-C}) \exp\left( C \int_{t_*}^T \| \sqrt{\n} \dot{u} \|_{L^2}^2 dt \right) \\
& \le C (1 + \tau^{-C}) \exp\left( C \int_{\frac{\tau}{2}}^T \| \sqrt{\n} \dot{u} \|_{L^2}^2 dt \right) \\
& \le C (1 + \tau^{-C}) e^{C\tau^{-C}},
\ea\ee
which together with (\ref{ngb121}) and the fact that $t_*<\tau$ gives
\be\ba\la{ngb126}
\sup_{\tau \le t \le T} \int \n |\dot u|^2 dx
\le \sup_{t_* \le t \le T} (2 f(t)) + C + C \tau^{-C}
\le C (1 + \tau^{-C}) e^{C\tau^{-C}}.
\ea\ee
It follows from (\ref{emdu}), (\ref{gb116a}), (\ref{nsb1}), (\ref{nsb02}), (\ref{nsb006}), (\ref{ngb125}), and $t_*<\tau$ that
\be\ba\la{ngb127}
\int_\tau^{T} \| \na\dot u\|^2_{L^2} dt
& \le C \int_{t_*}^T \left( \| \div \du \|^2_{L^2} + \| \curl \du \|^2_{L^2} 
+ \| \na u \|^4_{L^4} \right) dt \\
& \le C \int_{t_*}^T \left( \frac{1}{\nu^2} \| \dot{G} \|^2_{L^2} + \| \curl \du \|^2_{L^2} 
+ \| \na u \|^4_{L^4} + \| \na u \|^2_{L^2} \right) dt \\
& \le C (1 + \tau^{-C}) e^{C\tau^{-C}}.
\ea\ee
This combination of (\ref{ngb126}) and (\ref{ngb127}) yields (\ref{gb01a}).

It remains to prove (\ref{gb002}).
First, we deduce from (\ref{gb001}) and (\ref{nsb06}) that
there exists a positive constant $C$ depending only on $\mu,\ \ga,\ E_0,\ \|\n_0\|_{L^\infty},\ \| \na u_0 \|_{L^2},\ A$ and $\OM$ such that
\be\la{gb121}\ba
\sup_{0 \le t \le T} B^2_1
+ \int_0^T \left( B^2_1 + \| \sqrt{\n} \du \|^2_{L^2} \right) dt \le C.
\ea\ee
Multiplying (\ref{gb116b}) by $\si$ and using (\ref{nsb04}), (\ref{gb11006}), (\ref{gb121}), and Poincar\'e's inequality, we obtain
\be\la{gb122}\ba
& \frac{d}{dt} \left( \si f(t) \right)
+\frac{1}{4} \si \left( \frac{1}{\nu} \| \dot{G} \|^2_{L^2} + \mu \| \curl \dot{u}\|^2_{L^2} \right) \\
& \le \si'(t) f(t) + C \si \| \sqrt{\n} \du \|^4_{L^2} + C B^2_1 + C \| \sqrt{\n} \du \|^2_{L^2}.
\ea\ee
Moreover, by using (\ref{gb11006}), (\ref{gb118}) and (\ref{gb121}), we obtain
\be\la{gb123}\ba
\frac{1}{2} \| \sqrt{\n} \du \|^2_{L^2} \le f(t) + C,
\ea\ee
which together with (\ref{gb122}) gives
\be\la{gb124}\ba
& \frac{d}{dt} \left( \si f(t) \right)
+\frac{1}{4} \si \left( \frac{1}{\nu} \| \dot{G} \|^2_{L^2} + \mu \| \curl \dot{u}\|^2_{L^2} \right) \\
& \le C \si f(t) \| \sqrt{\n} \du \|^2_{L^2} + C B^2_1 + C \| \sqrt{\n} \du \|^2_{L^2}.
\ea\ee
The combination of (\ref{gb124}), (\ref{gb123}), (\ref{gb121}) and Gr\"onwall's inequality implies that
\be\la{gb125}\ba
& \sup_{0 \le t \le T} \si \int \rho|\dot{u}|^2dx + \int_0^T \si 
\left( \frac{1}{\nu} \| \dot{G} \|^2_{L^2} + \mu \| \curl \dot{u}\|^2_{L^2} \right) dt 
\le C.
\ea\ee
Finally, similar to (\ref{gb120}), after using (\ref{emdu2}), (\ref{nsb04}), (\ref{gb121}) and (\ref{gb125}), we derive
\be\la{gb126}\ba
\int_0^T \si \| \na \du \|^2_{L^2}
& \le C \int_0^T \si \left( \| \div \du \|^2_{L^2} + \| \curl \du \|^2_{L^2} 
+ \| \na u \|^4_{L^4} \right) dt \\
& \le C \int_0^T \si \left( \frac{1}{\nu^2} \| \dot{G} \|^2_{L^2} + \| \curl \du \|^2_{L^2} 
+ \| \na u \|^4_{L^4} + \| \na u \|^2_{L^2} \right) dt \\
& \le C,
\ea\ee
which together with (\ref{gb125}) yields (\ref{gb002}) and completes the proof of Lemma \ref{g1}.
\end{proof}

Adapting the arguments similar to those in \cite{LX3}, we can derive the following exponential decay estimate.
\begin{lemma}\la{bwle}
For any $s\in [1,\infty)$, there exist positive constants $C$, $K_0$, $\widetilde{\nu_0}$, where $C$ depends only on
$s$, $A$, $\OM$, $\ga$,\ $\mu$, $E_0$, $\|\n_0\|_{L^1\cap L^\infty}$, and $\| u_0\|_{H^1}$ such that if $\nu\geq\widetilde{\nu_0}$, it holds that
\be\la{be01}\ba
\sup_{1 \le t\le T} e^{\alpha_0 t} \left( \| \n-\ol{\n}\|_{L^s} + \| \na u \|_{L^s} + \| \sqrt{\n} \dot u \|_{L^2} \right)
\le C.
\ea\ee
where $\alpha_0 = \frac{K_0}{\nu}$.
\end{lemma}

\begin{lemma}\la{s22}
For any $2 < p < \infty$, there exists a positive constant $C$ depending only on 
$T$, $p$, $\ga$, $\mu$, $\nu$, $E_0$, $A$, $\OM$, $\| \rho_0 \|_{L^1 \cap L^\infty }$, and $\| \na u_0 \|_{L^2}$ such that
\be\la{s421} \ba
&\int_0^T \left(\|G\|_{L^\infty}+\|\na G\|_{L^p}+\|\o\|_{L^\infty}+\|\na \o\|_{L^p}+\| \rho \dot u\|_{L^p}\right)^{1+1 /p}dt \\& +\int_0^Tt\left(\|\na G\|_{L^p}^2+\|\na \o\|_{L^p}^2+\|   \dot u\|_{H^1}^2\right) dt \le C. \ea\ee
\end{lemma}
\begin{proof}
First, we deduce from (\ref{nsb06}) that
\be\la{s418}\ba
\sup_{0\le t\le T} \left( \| \n \|_{L^\infty} + \| u\|_{H^1} \right) + \int_0^T \left( \| \na u\|^2_{L^2}+ \| \sqrt{\n} \dot{u} \|^2_{L^2} \right)dt \le C.
\ea\ee

Then, multiplying (\ref{gb116}) by $\si$, choosing $\ep$ suitably small
and integrating over $(0,T)$, we use (\ref{gb118}), 
(\ref{s418}), (\ref{nsb04}) and (\ref{emdu2}) to derive
\be\la{s411}\ba
\sup_{0\le t\le T} \si \int\n|\dot u|^2dx+\int_0^{ T} \si \|\na\dot u\|^2_{L^2}dt\le C.
\ea\ee

In addition, H\"{o}lder's inequality and (\ref{pt1}) give
\be\la{s422} \ba
\| \rho \dot u\|_{L^p} & \le
C\| \rho \dot u\|_{L^2}^{2(p-1)/(p^2-2)}\|\dot u\|_{L^{p^2}}^{p(p-2)/(p^2-2)}\\ & \le
C\| \rho \dot u\|_{L^2}^{2(p-1)/(p^2-2)}\|\dot u\|_{H^1}^{p(p-2)/(p^2-2)}\\ & \le
C\| \rho^{1/2}  \dot u\|_{L^2} +C\| \rho \dot u\|_{L^2}^{2(p-1)/(p^2-2)}\|\na \dot u\|_{L^2}^{p(p-2)/(p^2-2)},
\ea\ee
which together with (\ref{s411}), (\ref{s422}) and (\ref{pt1}) yields that
\be\la{s423}\ba
&\int_0^T \left(\| \rho \dot u\|^{1+1 /p}_{L^p}+t\|  \dot u\|^2_{H^1}\right) dt\\
&\le C+C \int_0^T\left( \| \rho^{1/2}  \dot u\|_{L^2}^2 +  t\|\na \dot u\|_{L^2}^2+ 
t^{-(p^3-p^2-2p)/(p^3-p^2-2p+2)} \right)dt \\ 
&\le C.
\ea\ee
Moreover, from (\ref{s418}) and Sobolev embedding, we obtain
\be\la{s424}\ba 
& \|\div u\|_{L^\infty}+\|\o\|_{L^\infty} +  \|G\|_{L^\infty} \\ 
& \le C +C \|\na G\|_{L^p} +C \|\na \o\|_{L^p}\\ &\le C +C \|\n\dot u\|_{L^p}, \ea\ee 
which along with (\ref{nsb44}), (\ref{nsb04}), (\ref{s423}) and (\ref{s424}) implies (\ref{s421})
and completes the proof of Lemma \ref{s22}.
\end{proof}

\begin{lemma}\la{s23}
There exists a positive constant $C$ depending only on 
$T,\ q,\ \ga,\  \mu,\ \nu$, $E_0$, $A$, $\OM$, $\| \rho_0 \|_{L^1 \cap W^{1,q}}$ and $\| \na u_0 \|_{L^2}$, such that
\be\la{s431} \ba
&\sup_{0\le t\le T}\left(\norm[W^{1,q}]{ \rho}+\| u\|_{H^1}+t^{1/2}\| \n^{1/2} u_t \|_{L^2}+t^{1/2} \| u\|_{H^2} + \| \n_t \|_{L^2} \right) \\
& \quad + \int_0^T \left( \|\nabla^2 u\|^{2}_{L^2}+\|\nabla^2 u\|^{(q+1)/q}_{L^q}+t \|\nabla^2 u\|_{L^q}^2+\| \n^{1/2} u_t \|^2_{L^2}+t\| u_t\|_{H^1}^2 \right) dt\le C. \ea\ee
\end{lemma}
\begin{proof}
First, differentiating $(\ref{ns})_1$ with respect to $x$ and multiplying the 
resulting equation by $q |\nabla\n|^{q-2} \na \n $, we derive
\be\la{s432}\ba
& (|\nabla\n|^q)_t + \text{div}(|\nabla\n|^qu)+ (q-1)|\nabla\n|^q\text{div}u  \\
&+ q|\nabla\n|^{q-2} \p_i\n \p_i u^j \p_j\n +
q\n|\nabla\n|^{q-2}\p_i\n  \p_i \text{div}u = 0.\ea \ee 
Integrating (\ref{s432}) over $\OM$ and applying H\"older's inequality yields that
\be\la{s433}\ba
\frac{d}{dt} \norm[L^q]{\nabla\n}  
&\le C \norm[L^{\infty}]{\nabla u}  \norm[L^q]{\nabla\n} +C\|\na^2 u\|_{L^q} \\ 
&\le C(1+\norm[L^{\infty}]{\nabla u} ) \norm[L^q]{\nabla\n}+C \|\n\dot u\|_{L^q}, \ea\ee 
where we have used the following estimate:
\be\la{s434}\ba
\|\na^2 u\|_{L^q} &\le C(\|\na \div u\|_{L^q}+\|\na \o\|_{L^q}) \\
&\le C (\| \nabla G\|_{L^q}+ \|\nabla P \|_{L^q})+ C\|\na \o\|_{L^q} \\
&\le C\|\n\dot u\|_{L^q} + C  \|\nabla \n \|_{L^q} , \ea \ee 
owing to  (\ref{nsb44}) and (\ref{dc1}).

Additionally, by using (\ref{gn12}) and (\ref{s418}), we obtain
\be\la{s435}\ba
& \|\div u\|_{L^\infty}+\|\o\|_{L^\infty} \\
& \le C +C \|\na G\|_{L^q}^{q/(2(q-1))} + C \|\na \o\|_{L^q}^{q/(2(q-1))} \\
& \le C +C \|\n\dot u\|_{L^q}^{q/(2(q-1))},
\ea\ee
which together with Lemma \ref{bkm} and (\ref{s434}) implies that
\be\la{s436}\ba   
\|\na u\|_{L^\infty }  
&\le C\left(\|{\rm div}u\|_{L^\infty }+ \|\o\|_{L^\infty } \right)\log(e+\|\na^2 u\|_{L^q}) +C\|\na u\|_{L^2} +C \\
&\le C\left(1+\|\n\dot u\|_{L^q}^{q/(2(q-1))}\right)\log(e+\|\rho \dot u\|_{L^q} +\|\na \rho\|_{L^q}) +C\\
&\le C\left(1+\|\n\dot u\|_{L^q} \right)\log(e+ \|\na \rho\|_{L^q}) . \ea\ee
We set
\bnn
f(t) \triangleq e + \|\na \rho\|_{L^q},\quad h(t) \triangleq 1 + \|\n\dot u\|_{L^q},
\enn
which along with (\ref{s433}) and (\ref{s436}) shows that
\be\la{s4360}\ba
f'(t)\le Ch(t)f(t)\log f(t), 
\ea\ee 
due to $f(t)\ge e$ and $h(t)\ge 1$.

Dividing (\ref{s4360}) by $f(t)$ yields
\be\la{s437}\ba
(\log f(t))'\le Ch(t)\log f(t).
\ea\ee
From (\ref{s421}), we have
\be\la{s438}\ba 
\int_0^T \| \rho \dot u\|_{L^q}^{1+1/q}  dt \le C.
\ea\ee
Consequently, we conclude from (\ref{s437}), (\ref{s438}) and Gr\"onwall's inequality that
\be\la{s439}\ba
\sup\limits_{0\le t\le T} \| \nabla \rho \|_{L^q} \le C,
\ea\ee
which together with (\ref{s411}), (\ref{s421}), (\ref{s434}) and (\ref{s438}) leads to
\be\la{s4310}\ba
\sup_{0\le t\le T} t^{1/2} \| \na^2 u\|_{L^2} + \int_0^T \left(\|\nabla^2 u\|^{2}_{L^2}
+ \| \nabla^2 u \|^{(q+1)/q}_{L^q}+t \|\nabla^2 u\|_{L^q}^2 \right) dt \le C.
\ea\ee

In addition, it follows from $(\ref{ns})_1$, (\ref{s418}) and (\ref{s439}) that
\bnn  
\| \n_t\|_{L^2}\le
C\|u\|_{L^{2q/(q-2)}}\|\nabla \n \|_{L^q}+C\|\n\|_{L^\infty} \|\nabla u\|_{L^2} \le C,\enn
which gives
\be\la{s4311}\ba
\sup_{0\le t\le T} \| \n_t\|_{L^2} \le C.
\ea\ee

Finally, by using (\ref{s418}) and H\"older's inequality, we derive
\be \la{s4312}\ba 
\int\rho|u_t|^2dx 
&\le \int\rho| \dot u |^2dx+\int \n |u\cdot\na u|^2dx \\
&\le \int\rho| \dot u |^2dx+C \| u \|_{L^4}^2 \| \na u\|_{L^4}^2 \\ 
&\le \int\rho| \dot u |^2dx+C\| \na^2 u\|_{L^2}^2,
\ea\ee
and
\be\la{s4313}\ba 
\|\nabla u_t\|_{L^2}^2 
&\le \| \nabla \dot u \|_{L^2}^2+ \| \nabla(u\cdot\nabla u)\|_{L^2}^2  \\ 
&\le \|\nabla \dot u\|_{L^2}^2+ \|u\|_{L^{2q/(q-2)}}^2\|\nabla^2u \|_{L^q}^2+ \| \nabla u \|_{L^4}^4 \\ 
&\le \|\nabla \dot u\|_{L^2}^2+C\|\nabla^2u \|_{L^q}^2+ \| \nabla u \|_{L^4}^4.
\ea\ee 
Therefore, we deduce from (\ref{pt1}), (\ref{s411}), (\ref{s4310}), (\ref{s4312}) and (\ref{s4313}) that
\be\la{s4314}\ba
\sup_{0\le t\le T} t^{1/2}\| \n^{1/2} u_t \|_{L^2} + \int_0^T \| \n^{1/2} u_t \|^2_{L^2}
+ t\|u_t\|_{H^1}^2 dt \le C.
\ea\ee

By combining (\ref{s418}), (\ref{s439}), (\ref{s4310}), (\ref{s4311}) and (\ref{s4314}),
we obtain (\ref{s431}) and the proof of Lemma \ref{s23} is completed.
\end{proof}

From now on, we assume that the initial data $(\n_0,u_0)$ satisfy (\ref{csol1}) and the compatibility condition (\ref{csol2}).
\begin{lemma}\la{c21}
There exists a positive constant $C$ depending only on  
$T,\ \ga,\ \mu,\ \nu,\ E_0$, $A$, $\OM$, $\| \rho_0 \|_{L^1 \cap L^\infty }$, 
$\| \na u_0 \|_{L^2}$ 
and $\|g\|_{L^2}$, such that
\be\la{c411} \ba
\sup_{0\le t\le T}
\int\n|\dot u|^2dx
+\int_0^{ T}\int|\na\dot u|^2dxdt\le C.
\ea\ee
\end{lemma}
\begin{proof}
Taking into account the compatibility condition (\ref{csol2}), we define
\be\ba\nonumber
\sqrt{\n} \dot{u}(x,t=0)=g(x).
\ea \ee
By integrating (\ref{gb116}) over $(0,T)$, choosing $\ep$ sufficiently small,
and applying (\ref{gb118}), (\ref{s418}), (\ref{nsb04}) and (\ref{emdu2}), we obtain (\ref{c411}).
\end{proof}

In order to extend the local classical solution to a global solution, we need the following higher-order estimates.
Since the proofs of these estimates are similar to those in \cite{CL,LX3}, we omit the proofs.
\begin{lemma}\la{c26}
There exists a positive constant C depending only on $T,\ \mu,\ \nu,\ \ga,\ E_0$, $\| \rho_0 \|_{L^1 \cap W^{2,q}}$, $\| P(\rho_0) \|_{W^{2,q}}$, $\| u_0 \|_{H^2}$ and $\|g\|_{L^2}$ such that
\be\la{c421}\ba
& \sup_{0\le t\le T}\left(\norm[W^{1,q}]{ \rho} + \| u\|_{H^2} +\| \n^{1/2} u_t\|_{L^2}+\| \n_t\|_{L^2} \right) \\
& \quad + \int_0^T \left( \|\nabla^2 u\|_{L^q}^2 + \| \na u_t\|^2_{L^2} \right) dt\le C,
\ea\ee
\be\la{c431}\ba 
& \sup_{t\in[0,T]} \left(\norm[H^2]{\rho} + \norm[H^2]{P(\rho)}+ \|\n_t\|_{H^1}+\|P_t\|_{H^1} \right) \\
&\quad + \int_0^T\left(\| \na^3 u\|^2_{L^2}+\|\n_{tt}\|_{L^2}^2 +\|P_{tt}\|_{L^2}^2\right)dt \le C,
\ea\ee
\be\la{c442}\ba
\sup\limits_{0\le t\le T} t^{1/2} \left(\| \na u_t\|_{L^2}+\| \na^3 u\|_{L^2} \right) 
+ \int_0^T t\left(\|\n^{1/2}u_{tt}\|^2_{L^2}+\| \na^2 u_t\|^2_{L^2}\right)dt
\le C,
\ea\ee
\be\la{c451} \ba
\sup_{0\leq t\leq T}\left(\|\nabla^2 \n\|_{L^q } +\|\nabla^2 P  \|_{L^q }\right) \leq C,
\ea\ee
\be\la{c461}\ba
\sup_{0\leq t\leq T} t\left(\|\n^{1/2}u_{tt}\|_{L^2}+  \|\na^3 u \|_{L^q} + \|\na^2
u_t \|_{L^2}  \right) +\int_{0}^T  t^2 \|\nabla u_{tt}\|_{L^2}^2 dt\leq C.
\ea\ee
\end{lemma}

\section{Proofs of Theorems \ref{th0}--\ref{th3}}
In this section, we are devoted to proving the main results. 
Since the a priori estimates in Proposition \ref{pl1} 
require the density to be strictly away from vacuum,
we first establish the global existence of 
the classical solution to problem (\ref{ns})--(\ref{i3}) 
in the absence of vacuum. 
For cases where the initial density allows vacuum, 
it is noteworthy that all the a priori estimates in Sections 3 
and 4 are independent of the lower bound of the initial density.
Therefore, we approximate the initial density to ensure
that it remains strictly positive.
Finally, by applying standard compactness arguments,
we can prove the global existence.
\begin{proposition}\la{51}
Assume $\left( \n _0 ,u_0  \right)$ satisfy that for some $q>2$
\be\la{y521}\ba
\n _0 \in W^{2,q}, \quad \inf\limits_{x\in\OM}\n_0(x) >0, \quad u_0 \in H^2 \cap \tilde{H}^1,
\ea \ee 
and the compatibility condition (\ref{csol2}). 
Then when $\nu \ge \nu_1$, the problem \eqref{ns}--\eqref{i3} admits a unique classical solution $(\n,u)$ in $\OM \times (0,\infty)$ 
satisfying (\ref{csol4}) and
\be\la{y52121}\ba
\inf\limits_{(x,t)\in\OM \times (0,\infty)}\n(x,t) \ge C_0 >0.
\ea\ee
\end{proposition}
\begin{proof}
By the local existence result Lemma \ref{lct}, 
there exists a $T_*>0$	such that the problem (\ref{ns})--(\ref{i3}) 
has a unique classical solution $(\n,u)$ on $\OM \times (0,T_*]$. 
Next, we use a priori estimates Proposition \ref{pl1} and          
Lemma \ref{c26}, to extend the local classical solution$(\n,u)$ to all time.

Firstly, since $\n \in C \left([0,T_*];W^{2,q} \right)$, 
and $\n_0$ satisfies 
\bnn\ba 0<\inf\limits_{x\in\OM}\n_0(x) \le \n_0 \le \|\n_0\|_{L^\infty},\ea\enn
there exists a $T_1 \in (0,T_*] $ such that (\ref{p11}) holds for $T=T_1$.

Next, we introduce the following notation:
\be \la{y522}\ba	
T^* \triangleq \sup\{T\,|\,(\ref{p11}) \  \text{holds} \}.
\ea \ee	
Obviously, $T^*\ge T_1>0$. 
Furthermore, for any $0<\tau<T\leq T^*$
with $T$ finite, we derive from Lemma \ref{c26} that
\be \la{y523}\ba 
u \in C \left([\tau ,T];C^2(\ol{\OM}) \right) ,\quad
u_t \in C\left([\tau ,T];C(\ol{\OM}) \right),\ea \ee 
where we have used the standard embedding
$$L^\infty(\tau ,T;W^{3,q})\cap H^1(\tau ,T;H^2)\hookrightarrow
C\left([\tau ,T];C^2(\ol{\OM}) \right),  $$
and	
$$L^\infty(\tau ,T;H^2)\cap H^1(\tau ,T;L^2)\hookrightarrow
C\left([\tau ,T]; C(\ol{\OM}) \right).  $$	
Moreover, it follows from $(\ref{ns})_1$, Lemma \ref{c26}
as well as the standard arguments in \cite{L1} that
\be \la{y524}\ba 
\n \in C \left([0,T];W^{2,q} \right).\ea \ee 	
By combining (\ref{y523}) with (\ref{y524}), we derive
\be\la{y526}\ba 
\n^{1/2} \dot u \in C([\tau,T];L^2).\ea \ee
Finally, we claim that
\be\la{y527}\ba T^*=\infty .\ea \ee
Assume, for the sake of contradiction, that  $T^*<\infty$.
Then by Proposition \ref{pl1}, we obtain that (\ref{p11}) holds for
$T=T^*$.
We deduce from Lemma \ref{c26}, (\ref{y523}), (\ref{y524}), (\ref{y526})
that $(\n(x,T^*),u(x,T^*))$ satisfies (\ref{csol1}) and (\ref{csol2}), 
where $g(x)\triangleq \left( \n^{1/2} \dot u \right)(x, T^*),\,\,x\in \OM$. 
By combining $(\ref{ns})_1$ with (\ref{c421}) and using standard calculations, we conclude that
\be\ba\nonumber
\n(x,T^*) \ge \inf\limits_{x\in\OM}\n_0(x)\exp\left\{-\int_0^{T^*}\|\div u\|_{L^\infty}dt\right\} > 0.\ea \ee
Thus, Lemma \ref{lct} implies that there exists some $T^{**}>T^*$, such that
(\ref{p11}) holds for $T=T^{**}$, which contradicts (\ref{y522}),
and hence (\ref{y527}) holds. 
Finally, Lemma \ref{c26} shows that $(\rho,u)$
is in fact the classical solution defined on $\OM \times(0,T]$ for any
$0<T<T^*=\infty$. 
Furthermore, we deduce from (\ref{c431}) and (\ref{y521}) that (\ref{y52121}), the proof of Proposition \ref{51} is finished.
\end{proof}

\noindent\textbf{Proof of Theorem \ref{th2}}.
Let $(\n_0,u_0)$ be the initial data in Theorem \ref{th2}, satisfying (\ref{csol1}).
For any $\de \in (0,1)$, we set
\be\la{czbj1}\ba
\n_0^{\de} \triangleq \left( \n_0^\ga + \de^\ga \right)^{\frac{1}{\ga}},
\ea\ee
which together with (\ref{csol1}) gives
\be\ba\nonumber
0<\de \le \n_0^\de \le \| \n_0 \|_{L^\infty}+1,
\ea\ee
and
\be\ba\nonumber
\lim\limits_{\de\rightarrow 0} \|\n_0^\de-\n_0\|_{ W^{2,q}} = 0.
\ea\ee
In addition, we define
\be\la{czbj3}\ba
g^\de_2 := ( \n^\de_0 )^{-1/2} \left( - \mu \Delta u_0 - (\mu + \lm )\nabla \div u_0 
+ \nabla P(\n^\de_0) \right),
\ea\ee
The compatibility condition (\ref{csol2}) immediately implies
\be\ba\nonumber
g^\de_2 = ( \n^\de_0 )^{-1/2} ( \n_0 )^{1/2} g,
\ea\ee
which along with (\ref{czbj1}) yields
\be\ba\nonumber
\| g^\de_2 \|_{L^2} \le \| g \|_{L^2}.
\ea\ee

According to Proposition \ref{51}, we conclude that the problem (\ref{ns})--(\ref{i3}),
where $(\n_0,u_0)$ are replaced by $(\n_0^\de,u_0^\de )$ 
and the compatibility condition (\ref{csol2}) is replaced by (\ref{czbj3}),
admits a unique global classical solution  $(\n^\de,u^\de)$ satisfying Lemmas \ref{c21}
and \ref{c26}, with all constants $C$ independent of $\de$.
By letting $\de\rightarrow 0$
and using standard arguments (see \cite{HL,LZZ}), we obtain that the problem  (\ref{ns})--(\ref{i3}) 
has a global classical solution $(\n,u)$ satisfying (\ref{ssol4}).
Moreover, by using (\ref{be01}), we derive $(\n,u)$ satisfies the estimate (\ref{wsol4}).
The proof of uniqueness of $(\n,u)$ satisfying (\ref{ssol4}) is similar to \cite{Ge},
and hence we complete the proof of Theorem \ref{th2}.

\noindent\textbf{Proofs of Theorems \ref{th0} and \ref{th1}}.
By employing standard compactness arguments in \cite{F,L2},
Theorem \ref{th0} and Theorem \ref{th1} can be proven similarly to Theorem \ref{th2},
and hence their proofs are omitted.

\noindent\textbf{Proof of Theorem \ref{th01}}.
For any $\nu \ge \nu_1$, we deduce from (\ref{wsol2}), (\ref{nsb1}), (\ref{nsb006}) and Poincar\'e's inequality that $\{\n^{\nu}\}_\nu$ is bounded in $L^\infty(\OM \times (0,\infty))$
and that for any $0<\tau<T<\infty$, $\{u^{\nu}\}_\nu$ is bounded in $L^2(0,\infty;H^1) \cap L^\infty(\tau,\infty;H^1)$.
Moreover, by (\ref{pt1}) and H\"older's inequality, we have
\be\ba\nonumber
\| u^{\nu}_t \|_{L^2}
& \le C \left( \| \dot{u^{\nu}} \|_{L^2} + \| u^{\nu} \cdot \na u^{\nu} \|_{L^2} \right) \\
& \le C \left( \| \sqrt{\n^{\nu}} \dot{u^{\nu}} \|_{L^2} + \| \na \dot{u^{\nu}} \|_{L^2}
+ \| u^{\nu} \|_{L^4} \| \na u^{\nu} \|_{L^4} \right),
\ea\ee
which together with (\ref{gb01a}), (\ref{nsb02}), (\ref{nsb006}), (\ref{nsb04}) and Poincar\'e's inequality yields that for any $\tau>0$, $\{u^{\nu}\}_\nu$ is bounded in $H^1(\tau,\infty;L^2)$.

Therefore, we may assume that there exist a sequence $\nu_n \to \infty$ and functions $\n$ and $u$ such that, setting
\be\nonumber\ba
\n^n \triangleq \n^{\nu_n}, \quad u^n \triangleq u^{\nu_n},
\ea\ee
we have, for any $1 \le p < \infty$ and $0<\tau<T<\infty$,
\be\la{ins1}\ba
\begin{cases}
\n^n \rightharpoonup \n  \mbox{ weakly* in } L^\infty(\om \times (0,\infty)),\\
u^n \rightharpoonup u  \mbox{ weakly* in } \ L^\infty(\tau,T;H^1)\cap L^2(0,T;H^1), \\
u^n \to u  \mbox{ strongly  in } \ L^\infty(\tau,T;L^p).
\end{cases}
\ea\ee

Define
\be\nonumber\ba
G^n \triangleq \nu_n \div u^n - (P^n-\ol{P^n}), \quad \o^n \triangleq \na^\bot \cdot u^n,
\ea\ee
where $P^n = P(\n^n)$.
Since the pair $(\n^n,u^n)$ satisfies
\be\ba\la{ins2}
\begin{cases}
(\n^n)_t+\div(\n^n u^n)=0, \\
(\n^n u^n)_t+\div(\n^n u^n\otimes u^n) -\na G^n -\mu \na^{\bot} \o^n = 0,
\end{cases}
\ea\ee
we deduce that $G^n$ and $\o^n$ satisfy the following elliptic equations, respectively:
\be\ba\nonumber
\begin{cases}
\Delta G^n=\div \left(\rho^n \dot{u}^n \right) &\mathrm{in}\, \,  \OM, \\
\frac {\p G^n}{\p n}=\left(\rho^n \dot{u}^n-\mu \na^{\bot} \o^n \right) \cdot n &\mathrm{on}\, \,  \p \OM.
\end{cases}
\ea\ee
and
\be\ba\nonumber
\begin{cases}
\Delta \o^n=\na^\bot \cdot \left(\rho^n \dot{u}^n \right) &\mathrm{in}\, \,  \OM, \\
\o^n=-Au^n \cdot n^\bot &\mathrm{on}\, \,  \p \OM.
\end{cases}
\ea\ee
According to the standard $L^2$ estimate and Poincar\'e's inequality, we have
\be\la{ins3}\ba
\| G^n \|_{H^1} + \| \o^n \|_{H^1}
\le C \left( \| \sqrt{\n^n} \dot{u}^n \|_{L^2} + \| \na u^n \|_{L^2} \right).
\ea\ee

It follows from (\ref{ins3}), (\ref{nsb1}), and (\ref{nsb006}) that $\{G^n\}_n$ and $\{\o^n\}_n$ are bounded in $L^2(\tau,\infty;H^1)$.
Thus, by a diagonal argument, after extracting a further subsequence if necessary, there exists $\pi \in L^2_{ \text{loc} }(0,\infty;H^1)$ such that for any $\tau>0$,
\be\ba\la{ins4}
G^n \rightharpoonup - \pi, \quad \o^n \rightharpoonup \o \quad  \mbox{ weakly in } L^2(\tau,\infty;H^1),
\ea\ee
where $\o = \na^\bot \cdot u$.

Using (\ref{ins1}) and (\ref{ins4}), we may pass to the limit in (\ref{ins3}) as $n \to \infty$ in the sense of distributions and obtain that $(\n,u)$ satisfies
\be\ba\nonumber
\begin{cases}
\n_t+\div(\n u)=0,\\
(\n u)_t+\div(\n u\otimes u) -\mu \na^{\bot} \o + \na \pi =0.
\end{cases} 
\ea\ee
Moreover, (\ref{nsb1}) and (\ref{nsb006}) imply that for any $0<\tau<\infty$,
\be\la{ins12}\ba
\div u^\nu \to 0  \  \mbox{ strongly in } L^2(\OM \times (0,\infty)) \cap L^\infty(\tau,\infty;L^2),
\ea\ee
which gives (\ref{isol3}).
Then, from (\ref{ins1}) and (\ref{ins12}), we deduce that $\div u = 0$.
Since $\Delta u = \na \div u + \na^{\bot} \o $, it follows that $\na^{\bot} w=\Delta u$.
Hence, $(\n,u)$ is a global weak solution of (\ref{isol1}) and satisfies (\ref{isol2}).

We next consider the case of divergence-free initial velocity.
Assume, in addition, that the initial data $(\n_0,u_0)$ satisfy (\ref{insc1}).
Then, for any $\nu \ge \nu_1$, the a priori estimates (\ref{wsol2}), (\ref{nsb1}), (\ref{nsb06}), (\ref{gb002}), and Poincar\'e's inequality, together with $\div u_0 = 0$, imply that, for any $0<T<\infty$, $\{\n^{\nu}\}_\nu$ is bounded in $L^\infty(\om \times (0,\infty))$,
$\{u^{\nu}\}_\nu$ is bounded in $L^2(0,\infty;H^1) \cap L^\infty(0,\infty;H^1)$,
$\{\sqrt{\n^\nu} \du^{\nu}\}_\nu$ is bounded in $L^2(\om \times (0,\infty))$,
$\{\sqrt{t} \sqrt{\n^\nu} \dot{u^{\nu}}\}_\nu$ is bounded in $L^\infty(0,T;L^2)$,
and $\{\sqrt{t} \na \dot{u^{\nu}}\}_\nu$ is bounded in $L^2(\om \times (0,T))$.
Arguing as above, we obtain a subsequence of $(\n^{\nu},u^{\nu})$ converging to a global solution of (\ref{isol1}) satisfying (\ref{insc3}).

Next, we prove uniqueness for solutions of (\ref{isol1}) in the regularity class (\ref{insc3}).
Let $(\n_1,u_1,\pi_1)$ and $(\n_2,u_2,\pi_2)$ be two such solutions with the same initial data $(\n_0,u_0)$.
Set
\be\nonumber\ba
\de \n \triangleq \n_1-\n_2, \quad \de u \triangleq u_1-u_2, \quad \de \pi \triangleq \pi_1-\pi_2.
\ea\ee
A direct calculation shows that $(\de \n,\de u,\de \pi)$ satisfies
\be\la{rjwyx1}\ba
\begin{cases}
\de\n_t + u_2 \cdot \na \de\n + \de u \cdot \na\n_1 = 0, \\
\n_1 (\de u_t + u_1 \cdot \na \de u) - \mu \Delta \de u + \na \de \pi = -\de \n \dot{u}_2 - \n_1 \de u \cdot \na u_2, \\
\div \de u=0, \\
\de u \cdot n = 0, \quad \curl \de u=-A \de u \cdot n^\bot \quad \textnormal{on}\,\,\, \partial\Omega, \\
(\de \n,\de u)|_{t=0} = 0,
\end{cases} 
\ea\ee
where
\be\nonumber\ba
\dot{u}_2 \triangleq (u_2)_t + u_2 \cdot \na u_2.
\ea\ee

We first prove the uniqueness on a sufficiently short time interval.
Fix $0<T<1$.
Consider the Neumann problem
\be\la{rjwyx2}\ba
\begin{cases}
-\Delta \varphi = \de \n, \ & \textnormal{ in } \Omega, \\
\na \varphi \cdot n = 0, \  & \textnormal{ on } \partial\Omega, \\
\int_{\OM} \varphi \ dx = 0.
\end{cases}
\ea\ee
Since
\be\nonumber\ba
\int_{\OM} \de \n \ dx = 0,
\ea\ee
standard elliptic theory (see \cite{EL,NS}) guarantees the solvability of this system.

Multiplying $(\ref{rjwyx2})_1$ by $\varphi$, integrating by parts over $\OM$, and using the boundary condition $\na \varphi \cdot n = 0$ on $\p \OM$ and Poincar\'e's inequality, we obtain
\be\ba\la{rjwyx3}
\int |\na \varphi|^2 dx = \int \de \n \varphi dx
\le C \| \de \n \|_{L^2} \| \na \varphi \|_{L^2},
\ea\ee
which implies
\be\la{rjwyx4}\ba
\| \na \varphi \|_{L^2} \le C \| \de \n \|_{L^2}.
\ea\ee
Multiplying $(\ref{rjwyx1})_1$ by $\varphi$ and integrating by parts over $\OM$, we derive
\be\la{rjwyx5}\ba
\frac{1}{2} \frac{d}{dt} \| \na \varphi \|^2_{L^2}
& = \int u_2 \cdot \na \Delta \varphi \varphi dx - \int \de u \cdot \na \n_1 \varphi dx \\
& = \int \p_j u^i_2 \p_j \varphi \p_i \varphi dx + \int \de u \cdot \na \varphi \n_1 dx \\
& \le C \| \na u_2 \|_{L^\infty} \| \na \varphi \|^2_{L^2}
+ C \| \sqrt{\n_1} \de u \|_{L^2} \| \na \varphi \|_{L^2},
\ea\ee
which shows
\be\la{rjwyx6}\ba
\frac{d}{dt} \| \na \varphi \|_{L^2}
\le C \| \na u_2 \|_{L^\infty} \| \na \varphi \|_{L^2} + C \| \sqrt{\n_1} \de u \|_{L^2}.
\ea\ee
Integrating (\ref{rjwyx6}) over $(0,t)$ yields
\be\la{rjwyx7}\ba
\| \na \varphi(t,\cdot) \|_{L^2}
\le C \int_0^t \| \na u_2(s,\cdot) \|_{L^\infty} \| \na \varphi(s,\cdot) \|_{L^2} ds
+ C \int_0^t \| \sqrt{\n_1} \de u(s,\cdot) \|_{L^2} ds,
\ea\ee
where we have used $\na \varphi(0)=0$, which follows from (\ref{rjwyx4}) and $\de \n(0)=0$.

We set
\be\la{rjwyx8}\ba
D(t) \triangleq \sup_{0 < s \le t} s^{-\frac{1}{2}} \| \na \varphi(s,\cdot) \|_{L^2}.
\ea\ee
Combining this with (\ref{rjwyx7}) and H\"older's inequality leads to
\be\la{rjwyx9}\ba
D(t) \le C \int_0^t \| \na u_2(s,\cdot) \|_{L^\infty} D(s) ds + C \| \sqrt{\n_1} \de u \|_{L^2(\OM \times (0,t))}.
\ea\ee
The regularity in (\ref{insc3}), together with (\ref{gn11}) and H\"older's inequality, implies that
\be\la{rjwyx9a}\ba
\na u_2 \in L^1(0,T;L^\infty).
\ea\ee
By virtue of (\ref{rjwyx9}), (\ref{rjwyx9a}), and Gr\"onwall's inequality, it holds that
\be\la{rjwyx10}\ba
D(t) \le C \| \sqrt{\n_1} \de u \|_{L^2(\OM \times (0,t))}.
\ea\ee

We next estimate $\| \sqrt{\n_1} \de u \|_{L^2(\OM \times (0,t))}$.
To this end, consider the backward Stokes system
\be\la{rjwyx11}\ba
\begin{cases}
\n_1 (v_t + u_1 \cdot \na v) + \na \Pi + \mu \Delta v = \n_1 \de u, \\
\div v = 0, \\
v \cdot n =0, \quad \curl v = -A v \cdot n^\bot \quad \  \text{on} \  \p \OM, \\
v|_{t=T}=0.
\end{cases}
\ea\ee

If $\n_1$ and $u_1$ are regular with $\n_1$ bounded away from zero, the solvability of (\ref{rjwyx11}) follows from standard parabolic theory.
In the present setting, where vacuum may occur, one may replace $\n_1$ by $\n_1+\ep$, solve the corresponding nondegenerate backward problem, and then let $\ep \to 0$.
The a priori estimates derived below are uniform with respect to $\ep$, and hence standard compactness arguments, as in \cite{L1}, yield a solution of (\ref{rjwyx11}).

We now derive a priori estimates for system (\ref{rjwyx11}).
Multiplying $(\ref{rjwyx11})_1$ by $v$ and integrating by parts over $\OM \times (t,T)$, we obtain after using Young's inequality that
\be\la{rjwyx12}\ba
& \frac{1}{2} \| \sqrt{\n_1} v(t,\cdot) \|^2_{L^2}
+ \mu \int_t^T \int_{\OM} (\curl v)^2 dx d\tau
+ \mu \int_t^T \int_{\p \OM} A |v|^2 ds d\tau \\
& = - \int_t^T \int_{\OM} \n_1 \de u \cdot v dx d\tau
\le \frac{T}{4} \sup_{0 \le t \le T} \| \sqrt{\n_1} v \|^2_{L^2}
+ C \| \sqrt{\n_1} \de u \|^2_{L^2(\OM \times (0,T))},
\ea\ee
which together with (\ref{dc1}) and the fact that $T<1$ implies that
\be\la{rjwyx13}\ba
& \sup_{0 \le t \le T} \| \sqrt{\n_1} v \|^2_{L^2}
+ \int_0^T \|\na v\|_{L^2}^2 dt
\le C \| \sqrt{\n_1} \de u \|^2_{L^2(\OM \times (0,T))}.
\ea\ee

Define
\be\nonumber\ba
\dot{v} \triangleq v_t + u_1 \cdot \na v.
\ea\ee
Multiplying $(\ref{rjwyx11})_1$ by $\dot{v}$ and integrating by parts over $\OM$ yield
\be\la{rjwyx14}\ba
\| \sqrt{\n_1} \dot{v} \|^2_{L^2}
= - \int \na \Pi \cdot (v_t + u_1 \cdot \na v) dx
- \mu \int \Delta v \cdot (v_t + u_1 \cdot \na v) dx
+ \int \n_1 \de u \cdot \dot{v} dx.
\ea\ee

Next, we estimate each term on the right-hand side of (\ref{rjwyx14}).
Integration by parts combined with (\ref{gn11}) and Young's inequality implies
\be\la{rjwyx15}\ba
\left| \int \na \Pi \cdot (v_t + u_1 \cdot \na v) dx \right|
& = \left| \int u_1 \cdot \na v \cdot \na \Pi dx \right| \\
& \le \| \na \Pi \|_{L^2} \| \na v \|_{L^4} \| u_1 \|_{L^4} \\
& \le C \| \na \Pi \|_{L^2} \| \na v \|^{1/2}_{L^2} \| \na v \|^{1/2}_{H^1} \| \na u_1 \|_{L^2} \\
& \le \ep \left( \| \na \Pi \|_{L^2}^2 + \| \na^2 v \|_{L^2}^2 \right)
+ C(\ep) \| \na v \|^2_{L^2},
\ea\ee

\be\la{rjwyx16}\ba
\left| \mu \int u_1 \cdot \na v \cdot \Delta v dx \right|
& \le C \| \Delta v \|_{L^2} \| \na v \|_{L^4} \| u_1 \|_{L^4} \\
& \le C \| \Delta v \|_{L^2} \| \na v \|^{1/2}_{L^2} \| \na v \|^{1/2}_{H^1} \| \na u_1 \|_{L^2} \\
& \le \ep \| \na^2 v \|_{L^2}^2 + C(\ep) \| \na v \|^2_{L^2},
\ea\ee
and
\be\la{rjwyx17}\ba
\left| \int \n_1 \de u \cdot \dot{v} dx \right|
& \le \frac{1}{2} \| \sqrt{\n_1} \dot{v} \|^2_{L^2}
+ \frac{1}{2} \| \sqrt{\n_1} \de u \|^2_{L^2}.
\ea\ee
Integration by parts also gives
\be\la{rjwyx18}\ba
\int \Delta v \cdot v_t dx = \int \na^\bot \curl v \cdot v_t dx
= -\frac{1}{2} \frac{d}{dt} \left( \int (\curl v)^2 dx + \int_{\p \OM} A |v|^2 ds \right).
\ea\ee
Substituting (\ref{rjwyx15})--(\ref{rjwyx18}) into (\ref{rjwyx14}) shows
\be\la{rjwyx19}\ba
& - \frac{d}{dt} \left( \mu \int (\curl v)^2 dx + \mu \int_{\p \OM} A |v|^2 ds \right) + \| \sqrt{\n_1} \dot{v} \|^2_{L^2} \\
& \le \| \sqrt{\n_1} \de u \|^2_{L^2} + 4 \ep \left( \| \na \Pi \|_{L^2}^2 + \| \na^2 v \|_{L^2}^2 \right)
+ C(\ep) \| \na v \|^2_{L^2}.
\ea\ee

Since $\curl v$ satisfies
\be\la{rjwyx20}\ba
\begin{cases}
\mu \Delta \curl v = \na ^\bot \cdot \left( \n_1 \de u - \n_1 \dot{v} \right) & \text{ in } \OM, \\
\curl v =-A v \cdot n^{\bot} & \text{ on } \p\OM,
\end{cases}
\ea\ee
standard energy estimates and Poincar\'e's inequality yield
\be\la{rjwyx21}\ba
\| \na \curl v \|_{L^2} & \le C \left( \| \n_1 \de u \|_{L^2} + \| \n_1 \dot{v} \|_{L^2} + \| \na v \|_{L^2} + \| v \|_{L^2} \right) \\
& \le C \left( \| \sqrt{\n_1} \de u \|_{L^2} + \| \sqrt{\n_1} \dot{v} \|_{L^2} + \| \na v \|_{L^2} \right).
\ea\ee
This, together with (\ref{dc1}), implies
\be\la{rjwyx22}\ba
\| \na^2 v \|_{L^2} \le C \| \curl v \|_{H^1}
\le C \left( \| \sqrt{\n_1} \de u \|_{L^2} + \| \sqrt{\n_1} \dot{v} \|_{L^2} + \| \na v \|_{L^2} \right).
\ea\ee
Combining this with $(\ref{rjwyx11})_1$, we also arrive at
\be\la{rjwyx23}\ba
\| \na \Pi \|_{L^2} & \le C \left( \| \na^2 v \|_{L^2} + \| \n_1 \de u \|_{L^2} + \| \n_1 \dot{v} \|_{L^2} \right) \\
& \le C \left( \| \sqrt{\n_1} \de u \|_{L^2} + \| \sqrt{\n_1} \dot{v} \|_{L^2} + \| \na v \|_{L^2} \right).
\ea\ee
Putting (\ref{rjwyx22}) and (\ref{rjwyx23}) into (\ref{rjwyx19}) and choosing $\ep$ suitably small lead to
\be\la{rjwyx24}\ba
& - \frac{d}{dt} \left( \mu \int (\curl v)^2 dx + \mu \int_{\p \OM} A |v|^2 ds \right) + \frac{1}{2}  \| \sqrt{\n_1} \dot{v} \|^2_{L^2}
\le C \| \sqrt{\n_1} \de u \|^2_{L^2} + C \| \na v \|^2_{L^2}.
\ea\ee
Integrating (\ref{rjwyx24}) backward in time, using $v(\cdot,T)=0$, and applying (\ref{rjwyx13}), (\ref{rjwyx22}), and (\ref{rjwyx23}), we deduce that
\be\la{rjwyx25}\ba
\sup_{0\le t \le T} \| \na v \|^2_{L^2}
+ \int_0^T \left( \| \sqrt{\n_1} \dot{v} \|^2_{L^2} + \| \na^2 v \|^2_{L^2} + \| \na \Pi \|^2_{L^2} \right) dt
\le C \| \sqrt{\n_1} \de u \|^2_{L^2(\OM \times (0,T))}.
\ea\ee

Next, we multiply $(\ref{rjwyx1})_2$ and $(\ref{rjwyx11})_1$ by $v$ and $\de u$, respectively, integrate by parts, and sum the results to obtain
\be\la{rjwyx26}\ba
\| \sqrt{\n_1} \de u \|^2_{L^2(\OM \times (0,T))}
= \int_0^T \int \left( \de \n \dot{u}_2 \cdot v + \n_1 \de u \na u_2 \cdot v  \right) dx dt.
\ea\ee
Integration by parts along with (\ref{insc3}), (\ref{rjwyx8}), and H\"older's inequality yields
\be\la{rjwyx27}\ba
\int_0^T \int \de \n \dot{u}_2 \cdot v dx dt
& = - \int_0^T \int \Delta \varphi \dot{u}_2 \cdot v dx dt \\
& = \int_0^T \int \na \varphi \cdot \na (\dot{u}_2 \cdot v) dx dt \\
& \le C \int_0^T \| \na \varphi \|_{L^2} \left( \| \na \dot{u}_2 \|_{L^2} \| v \|_{L^\infty} + \| \dot{u}_2 \|_{L^4} \| \na v \|_{L^4} \right) dt \\
& \le C D(T) \int_0^T \left( \| t^{\frac{1}{2}} \na \dot{u}_2 \|_{L^2} \| v \|_{L^\infty} + \| t^{\frac{1}{2}} \dot{u}_2 \|_{L^4} \| \na v \|_{L^4} \right) dt \\
& \le C D(T) \| \na v \|_{L^2(0,T;L^4)} \\
& \le C T^{\frac{1}{4}} \| \sqrt{\n_1} \de u \|^2_{L^2(\OM \times (0,T))},
\ea\ee
where in the last inequality we have used (\ref{rjwyx10}) and the following estimate:
\be\la{rjwyx28}\ba
\int_0^T \| \na v \|^2_{L^4} dt \le C \int_0^T \| \na v \|_{L^2} \| \na v \|_{H^1} dt
\le C T^{\frac{1}{2}} \| \sqrt{\n_1} \de u \|^2_{L^2(\OM \times (0,T))}.
\ea\ee
due to (\ref{gn11}), (\ref{rjwyx13}), and (\ref{rjwyx25}).

On the other hand, it follows from (\ref{insc3}), (\ref{rjwyx28}), and H\"older's inequality that
\be\la{rjwyx29}\ba
\int_0^T \int \n_1 \de u \na u_2 \cdot v dx dt
& \le C \| \sqrt{\n_1} \de u \|_{L^2(\OM \times (0,T))}
\| \na u_2 \|_{L^\infty(0,T;L^2)} \| v \|_{L^2(0,T;L^\infty)} \\
& \le C T^{\frac{1}{4}} \| \sqrt{\n_1} \de u \|^2_{L^2(\OM \times (0,T))}.
\ea\ee
Inserting (\ref{rjwyx27}) and (\ref{rjwyx29}) into (\ref{rjwyx26}) gives
\be\la{rjwyx30}\ba
\| \sqrt{\n_1} \de u \|^2_{L^2(\OM \times (0,T))}
\le C T^{\frac{1}{4}} \| \sqrt{\n_1} \de u \|^2_{L^2(\OM \times (0,T))}.
\ea\ee
Thus, if $0<T<1$ sufficiently small, then
\be\la{rjwyx30a}\ba
\sqrt{\n_1} \de u = 0 \quad \text{ a.e. in } \OM \times (0,T).
\ea\ee
Combining this with (\ref{rjwyx10}), (\ref{rjwyx8}), and $(\ref{rjwyx2})_1$, we arrive at
\be\la{rjwyx30b}\ba
\de \n = 0 \quad \text{ a.e. in } \OM \times (0,T).
\ea\ee
Hence, $(\ref{rjwyx1})_2$ reduces to
\be\la{rjwyx31}\ba
\n_1 (\de u_t + u_1 \cdot \na \de u) - \mu \Delta \de u + \na \pi = 0.
\ea\ee
Multiplying (\ref{rjwyx31}) by $\de u$ and integrating by parts over $\OM \times (0,t)$, we obtain
\be\la{rjwyx32}\ba
\frac{1}{2} \| \sqrt{\n_1} \de u(t,\cdot) \|^2_{L^2}
+ \mu \int_0^t \int (\curl \de u)^2 dx d\tau
+ \mu \int_0^t \int_{\p \OM} A |\de u|^2 ds d\tau = 0,
\ea\ee
which together with (\ref{dc1}), $(\ref{rjwyx1})_3$, and (\ref{rjwyx30a}) yields
\be\nonumber\ba
\de u = 0 \quad \text{ a.e. in } \OM \times (0,T).
\ea\ee
Finally, from $(\ref{rjwyx1})_2$, we conclude that if the normalization
\be\nonumber\ba
\int_{\om} \pi(x,t) dx = 0,
\ea\ee
is imposed for almost every $t$, then $\de \pi = 0$ as well.

Repeating the preceding argument on successive sufficiently short time intervals proves uniqueness on every finite time interval.
Thus, $(\n_1,u_1,\pi_1)=(\n_2,u_2,\pi_2)$, and uniqueness follows.

By uniqueness, every convergent subsequence of $(\n^\nu,u^\nu)$ has the same limit, and hence the whole family $(\n^\nu,u^\nu)$ converges to the unique solution $(\n,u)$ of (\ref{isol1}) satisfying (\ref{isol2}).
This completes the proof of Theorem \ref{th01}.

\noindent\textbf{Proof of Theorem \ref{th3}}.
The proof of Theorem \ref{th3} is similar to the proof of \cite[Theorem 1.2]{CL}, so we omit it here.

\end{document}